\documentclass[12pt]{amsart}
\usepackage{amssymb}

\newtheorem{theorem}{Theorem}[section]
\newtheorem{lemma}[theorem]{Lemma}

\newtheorem{corollary}[theorem]{Corollary}
\newtheorem{proposition}[theorem]{Proposition}

\newtheorem{lem-def}[theorem]{Lemma-Definition}

\DeclareRobustCommand\longtwoheadrightarrow
{\relbar\joinrel\twoheadrightarrow}

\newcommand{\hooklongrightarrow}{\lhook\joinrel\longrightarrow}
\renewenvironment{proof}{{\bfseries Proof.}}{\qed}
\newcommand{\R}{\mathbb R}
\newcommand{\I}{\mathbb I}

\newcommand{\N}{\mathbb N}
\newcommand{\Z}{\mathbb Z}
\newcommand{\Q}{\mathbb Q}
\newcommand{\C}{\mathbb C}

\newcommand{\T}{\mathbb T}

\def\op{\operatorname}

\def\aa{\mathcal{A}}

\def\al{\alpha}

\def\ars#1{\renewcommand\arraystretch{#1}}

\def\aut{\op{Aut}}

\def\bb{{\mathcal B}}

\def\bs{\vskip.5cm}
\def\be{\beta}

\def\cc{\mathcal{C}}
\def\ccb{\mathcal{C}_{\op{bad}}}
\def\cfa{\left(\ri\right)_{i\in A}}
\def\cfb{\left(\zj\right)_{j\in B}}

\def\com{{\op{com}}}

\def\d{\Delta}
\def\defn{\nn{\bf Definition. }}
\def\dep{\op{\mbox{\small depth}}}

\def\dm{\Delta_\mu}

\def\diso{\lower.4ex\hbox{$\downarrow$}\raise.4ex\hbox{\mbox{\scriptsize
$\wr$}}}

\def\dta{\delta}

\def\e{\medskip}

\def\ee{\mathcal{E}}

\def\ep#1{\exp(\Pi i#1)}
\def\ep{\epsilon}

\def\erel{e_{\op{rel}}}

\def\fin{{\op{fin}}}

\def\g{\Gamma}
\def\ga{\gamma}

\def\gb{\g_\bb}
\def\gbc{\g_{\op{bc}}}

\def\gc{\g_\cc}

\def\gen#1{\big\langle\, {#1} \,\big\rangle}

\def\gg{\mathcal{G}}
\def\gga{\gen{\g,\ga}}
\def\ggb{\gen{\g,\be}}

\def\ggm{\mathcal{G}_\mu}

\def\ggn{\mathcal{G}_\nu}

\def\gm{\g_\mu}

\def\gn{\g_\nu}
\def\gnbc{\g_{\op{nbc}}}

\def\gq{\g_\Q}
\def\gr{\operatorname{gr}}
\def\grh{\g_\rho}

\def\gsme{\g_{\op{sme}}}

\def\hk{\hookrightarrow}
\def\hq{\mathbf{H}(\gq)}
\def\hra{\hooklongrightarrow}

\def\imp{\ \Longrightarrow\ }

\def\inii{\op{Init}(I)}
\def\init{\op{init}}
\def\inm{\op{in}_\mu}
\def\inn{\op{in}}

\def\inu{\op{in}_\nu}
\def\irr{\op{Irr}}
\def\ism{\lower.3ex\hbox{\ars{.08}$\begin{array}{c}\,\to\\\mbox{\tiny $\sim\,$}\end{array}$}}
\def\iso{\ \lower.3ex\hbox{\ars{.08}$\begin{array}{c}\lra\\\mbox{\tiny $\sim\,$}\end{array}$}\ }

\def\ka{\kappa}
\def\kam{\kappa(\mu)}
\def\kb{\overline{K}}
\def\kh{K^h}
\def\khx{K^h[x]}
\def\km{k_\mu}

\def\kp{\op{KP}}

\def\kpi{\op{KP}_\infty}
\def\kpm{\op{KP}(\mu)}
\def\kpmz{\op{KP}_{\op{str}}(\mu)}
\def\kpn{\op{KP}(\nu)}

\def\kpr{\op{KP}(\rho)}

\def\ks{K^{\op{sep}}}

\def\kx{K[x]}

\def\la{\lambda}
\def\La{\Lambda}
\def\lc{\op{lc}}
\def\lci{\ll\cc_\infty}
\def\ldp{\op{\mbox{\small lim-depth}}}
\def\lfin{\ll_\fin}
\def\lg{l\raise.6ex\hbox to.2em{\hss.\hss}l}
\def\li{\ll_\infty}
\def\ll{\mathcal{L}}
\def\lra{\,\longrightarrow\,}
\def\lui{\mathcal{LU}_\infty}
\def\lx{\operatorname{lex}}

\def\mi{m_\infty}
\def\minf{\om_{-\infty}}

\def\mlv{Maclane--Vaqui\'e\ }

\def\mmu{\mid_\mu}
\def\mon{\op{Mon}(L/K)}

\def\nb{\bar{\nu}}

\def\nn{\noindent}
\def\nni{\mathcal{N}_\infty}

\def\om{\omega}
\def\oo{\mathcal{O}}
\def\orb{\hbox to  .3em{$\backslash$}\backslash}
\def\ord{\op{ord}}
\def\p{\mathfrak{p}}
\def\paph{\pset_\aa(\phi)}

\def\pcv{\op{Prin}}

\def\pma{\pset(\mu_\aa)}
\def\pmi{\pset(\minf)}
\def\prh{\pset(\rho)}

\def\pmph{\pset_\mu(\phi)}

\def\ppa{\mathcal{P}_{\al}}
\def\prim{\op{Prim}}

\def\pset{\mathcal{P}}

\def\rha{\rho_\aa}
\def\rhb{\rho_\bb}
\def\rhc{\rho_\cc}
\def\ri{\rho_i}
\def\rii{\R^\I_{\lx}}
\def\rj{\rho_j}
\def\rk{\rho_k}
\def\rl{\rho_\ell}
\def\rll{\R_\sme}
\def\rlx{\R^I_{\lx}}
\def\rlxs{\R^{I_S}_{\lx}}

\def\rnk{\op{rk}}

\def\sg{\sigma}

\def\sii{\ \Longleftrightarrow\ }

\def\sme{{\mbox{\tiny $\op{sme}$}}}
\def\smu{\sim_\mu}

\def\str{^{\str}}

\def\spv{\op{Spv}}

\def\stan{_{\op{tan}}}
\def\str{^{\op{str}}}
\def\supp{\op{supp}}

\def\sval{\operatorname{sv}}
\def\t{\theta}
\def\tan{\sim_{\op{tan}}}

\def\tdm{\op{TD}(\mu)}
\def\tgq{\ttt(\gq)}

\def\tind{\ttt^{\op{ind}}_\sme}
\def\tla{\ttt(\La)}
\def\tmn{\ty(\mu,\nu)}
\def\tmr{\ty(\mu,\rho)}

\def\tq{\mathcal{T}_\Q}
\def\ts{\mathcal{T}_\sme}
\def\ttt{\mathcal{T}}
\def\ty{\mathbf{t}}
\def\tz{\mathcal{T}^0}

\def\vb{\bar{v}}

\def\wb{\overline{w}}

\def\zi{\zeta_i}
\def\zj{\zeta_j}
\def\zk{\zeta_k}

\newcounter{cs}
\stepcounter{cs}
\newcommand{\casos}{\begin{itemize}}
\newcommand{\fcasos}{\end{itemize}\setcounter{cs}{1}}

\newfont{\tit}{cmr12 scaled \magstep3}

\title[Valuative trees]{Valuative trees over valued fields}

\makeatletter
\@namedef{subjclassname@2010}{%
  \textup{2010} Mathematics Subject Classification}
%, 12J10}

\author[Alberich]{Maria Alberich-Carrami$\tilde{\mbox{n}}$ana}
\address{Institut de Rob\`otica i Inform\`atica Industrial (IRI, CSIC-UPC), Institut de Mate\-mà\-tiques de la UPC-BarcelonaTech (IMTech) and Departament de Mate\-m\`a\-tiques, Universitat Polit\`ecnica de Cata\-lunya $\cdot$ BarcelonaTech, Av. Diagonal, 647, E-08028 Barcelona, Catalonia}
\email{Maria.Alberich@upc.edu}
\author[Gu\`ardia]{Jordi Gu\`ardia}
\address{Departament de Matem\`atiques, Escola Polit\`ecnica Superior d'Enginye\-ria de Vilanova i la Geltr\'u, Av. V\'\i ctor Balaguer s/n. E-08800 Vilanova i la Geltr\'u, Catalonia}
\email{jordi.guardia-rubies@upc.edu}
\author[Nart]{Enric Nart}
\author[Ro\'e]{Joaquim Ro\'e}
\address{Departament de Matem\`{a}tiques,         Universitat Aut\`{o}noma de Barcelona,         Edifici C, E-08193 Bellaterra, Barcelona, Catalonia}
\email{nart@mat.uab.cat,\quad jroe@mat.uab.cat}
\thanks{Partially supported by grants  PID2020-116542GB-I00 and PID2019-103849GB-I00 from the Spanish Research Agency, and grant 2017SGR-932 from Generalitat de Catalunya}
\date{\today}

\begin{document}
%MTM2016-75980-P, MTM2015-69135-P

%\pagestyle{empty}

%No assumption is made on the rank of the valuations.
\begin{abstract}
For an arbitrary valued field $(K,v)$ and a given extension $v(K^*)\hk \La$ of ordered groups, we analyze the structure of the tree formed by all $\La$-valued extensions of $v$ to  the polynomial ring $\kx$.  As an application, we find a model for the tree of all equivalence classes of valuations on $\kx$ (without fixing their value group), whose restriction to $K$ is equivalent to $v$.
\end{abstract}

\maketitle

%\begin{center}\sl Preliminary version \end{center}

%\tableofcontents

\section*{Introduction}
A valuation on  a commutative ring $A$ is a mapping $\mu\colon A\to \La\infty$,
where $\La$ is an ordered abelian group, satisfying the following conditions:

(0) \ $\mu(1)=0$, \ $\mu(0)=\infty$,

(1) \ $\mu(ab)=\mu(a)+\mu(b),\quad\forall\,a,b\in A$,

(2) \ $\mu(a+b)\ge\min\{\mu(a),\mu(b)\},\quad\forall\,a,b\in A$.\e

The \emph{support} of $\mu$ is the prime ideal $\p=\supp(\mu)=\mu^{-1}(\infty)\in\op{Spec}(A)$.
%The valuation $\mu$ induces a valuation $\bar{\mu}$ on the field of fractions of $A/\p$.
The \emph{value group} of $\mu$ is the subgroup $\gm\subset \La$ generated by $\mu\left(A\setminus\p\right)$.

%Denote the maximal ideal, valuation ring and residue class field of this valuation on $\kappa(\p)$, by $$\m_\mu\subset \oo_\mu\subset \kappa(\p),\qquad k_{\mu}=\oo_\mu/\m_\mu.$$
%Note that $\kappa(0)=K(x)$, while for $\p\ne0$ the field $\kappa(\p)$ is a simple finite extension of $K$.

%Thus, a valuation $\mu$ on $\kx$ determines a valuation on a simple field extension $L/K$, either algebraic or transcendental.\bs

Two valuations $\mu$, $\nu$ on $A$ are \emph{equivalent} if there is an isomorphism of ordered groups $\iota\colon \gm \ism\gn$ such that $\nu=\iota\circ \mu$.
In this case, we write $\mu\sim\nu$.
%fitting into a commutative diagram$$\ars{1.3}\begin{array}{ccc}\gm\infty&\stackrel{\varphi}\lra\ &\!\!\gn\infty\\\quad\ \mbox{\scriptsize$\mu$}&\nwarrow\ \nearrow&\!\!\!\mbox{\scriptsize$\nu$}\quad\\&A&\end{array}$$

The \emph{valuative spectrum} of $A$ is the set $\spv(A)$ of equivalence classes of valuations on $A$. We denote by $[\mu]\in\spv(A)$ the equivalence class of $\mu$.

Any ring homomorphism $A\to B$ induces a restriction of valuations which behaves well on equivalence classes and determines a mapping $\spv(B)\to\spv(A)$. %determined by the assignment:$$B\stackrel{\mu}\lra\La\infty\quad\longmapsto\quad A\lra B\stackrel{\mu}\lra\La\infty,$$

%\subsection*{Extensions of a valuation on $K$}

For any field $K$ we may consider the relative affine line $\spv(\kx)\to\spv(K)$.

Given any valuation $v$ on $K$, the fiber $\ttt_v$ of the equivalence class $[v]\in\spv(K)$ is called the \emph{valuative tree} over the valued field $(K,v)$.
$$
\ars{1.2}
\begin{array}{ccc}
\ttt_v&\hra &\spv(\kx)\\
\downarrow&&\downarrow\\
\mbox{$[v]$}&\hra &\spv(K)
\end{array}
$$

This terminology is  borrowed from Favre and Jonsson's book \cite{FJ}, where valuations of certain 2-dimensional local rings are studied.
In the case $\rnk(\g)=1$ and $K$ algebraically closed, the valuative tree admits a structure of a Berkovich space and has relevant analytical properties \cite{Bch}.

The main aim of this paper is to obtain a thorough description of the tree $\ttt_v$ for an arbitrary valued field $(K,v)$.
 
In the first part of the paper, composed of sections 1--5, we fix an extension $\g\hk\La$ of ordered abelian groups, and we describe the tree $\ttt=\ttt(\La)$ formed by all $\La$-valued extensions of $v$ to $\kx$.

Section 1 includes some background on key polynomials of valued fields. For any valuation $\mu$ on $\kx$, let $\kpm$ be the set of Maclane-Vaqui\'e key polynomials for $\mu$. If $\kpm\ne\emptyset$, then $\mu$ has a \emph{degree} and a \emph{singular value}, defined as
$\deg(\mu)=\deg(\phi)$, $\sval(\mu)=\mu(\phi)$,
for any  key polynomial $\phi\in\kpm$  of minimal degree.

Section 2 discusses tangent directions and the tangent space of $\ttt$. The leaves of $\ttt$ (maximal nodes) are characterized by the property $\kpm=\emptyset$. The set of tangent directions of an inner node $\mu\in\ttt$ is parametrized by the set $\kpm/\!\smu$ of $\mu$-equivalence classes of key polynomials.

Section 3 describes the set $\lfin(\ttt)$ of \emph{finite leaves} of the tree, determined by all valuations with non-trivial support. There is a bijection between $\lfin(\ttt)$ and  the set of monic irreducible polynomials in $\kh[x]$, where $\kh$ is a henselization of $K$. This result is just a reformulation of classical valuation-theoretic results.

Section 4 describes the set of \emph{infinite leaves} of $\ttt$, which are a kind of limit of certain totally ordered families of inner nodes of $\ttt$. This section contains a detailed analysis of \emph{limit augmentations} of valuations too. This concept was introduced by Vaqui\'e in his fundamental papers  \cite{Vaq0,Vaq} extending  to arbitrary valued fields the pioneering work of Maclane for discrete rank-one valued fields \cite{mcla}.

Limit augmentations are based on \emph{continuous families} of iterated augmentations. In the literature, we find different conditions imposed on these families, serving different purposes. We define a continuous family as a totally ordered family of valuations in $\ttt$, containing no maximal element, and having a stable degree. There is a natural equivalence relation between these families and we show, in Lemma \ref{specialCont}, that every equivalence class of continuous families contains a family satisfying all relevant conditions that are attributed to these families in the literature.

Section 5 gives a detailed description of $\ttt$. Section 5.1 reviews the fundamental result of Maclane-Vaqui\'e describing how to reach all nodes of $\ttt$ by a combination of ordinary augmentations, limit augmentations and stable limits. Every node $\mu\in\ttt$ may be linked to some degree-one node in $\ttt$ by an essentially unique \emph{Maclane-Vaqui\'e (MLV) chain}, supporting  data intrinsically associated to $\mu$ \cite{MLV}.  For instance, each node $\mu\in\ttt$ has a \emph{depth}, defined as the length of its MLV chain, which is either a natural number or infinity.  In Section 5.2, we show that these intrinsic data encode  arithmetic or geometric invariants of $\mu$, depending on the context in which the base valued field $(K,v)$ is considered.   Sections  5.3--5.5, describe the different kinds of paths we may find in $\ttt$. In Section 5.6, we prove that every two nodes of $\ttt$ have a greatest common lower node in $\ttt$ and relate our description of $\ttt$ to the notion of $\Lambda$-tree.  

In the second part of the paper, composed of sections 6--7, we find a concrete model for the valuative tree $\ttt_v$. 

For any valuation $\mu$ on $\kx$ extending $v$, the embedding $\g\hk\gm$ is a \emph{small extension} of ordered groups; that is, if $\g'\subset\gm$ is the relative divisible closure of $\g$  in $\gm$, then the quotient $\gm/\g'$ is a cyclic group  \cite[Thm. 1.5]{Kuhl}.

In \cite{csme}, a universal extension $\g\hk\rii$ of ordered groups is constructed, which contains all small extensions of $\g$ up to $\g$-isomorphism as ordered groups. On a certain subset $\rll\subset\rii$,
an equivalence relation $\sim_\sme$ is defined such that the quotient set $\rll/\!\sim_\sme$ parametrizes the quasi-cuts of the divisible hull of $\g$. Also, there is a canonical subset $\gsme\subset\rll$ which faithfully represents all $\sim_\sme$ classes.

In Section 6, we consider the subtree $\tz\subset \ttt(\rii)$ formed by all nodes $\mu$ such that $\gm\subset \rll$. Then, we characterize equivalence of valuations in $\tz$ as follows.  \e

\noindent{\bf Proposition 6.3. }{\it Let $\mu,\nu\in\tz$ be two inner nodes. Then, $\mu\sim\nu$ if and only if the following three conditions hold:

(a) \ The valuations $\mu$, $\nu$ admit a common key polynomial of minimal degree.

(b) \ For all \,$a\in\kx\,$ such that \,$\deg(a)<\deg(\mu)$, we have $\,\mu(a)=\nu(a)$.

(c) \ $\sval(\mu)\sim_\sme \sval(\nu)$.

In this case, we have $\kpm=\kpn$.}\e

In Section 7, we consider the subtree $\ts\subset\tz$ formed by all leaves of $\tz$, and all inner nodes $\mu$ such that $\sval(\mu)$ belongs to $\gsme$. Then, we use Proposition 6.3 to obtain our main theorem.\e

\noindent{\bf Theorem 7.1 }{\it The mapping $\mu\mapsto[\mu]$ induces a bijection between $\ts$ and $\ttt_v$.}\e

In the rest of the section, we discuss special features of the paths in $\ts$ and we show the existence of \emph{primitive nodes}, leading to a certain stratification of the tree by \emph{limit-depth}, which is the number of limit augmentations in the MLV chains.\e

The techniques of this paper have been applied in two different contexts \cite{AGNR,Rig}. Let $(\kh,v^h)$ be a henselization of $(K,v)$. In \cite{AGNR}, we use the primitive nodes of the valuative tree to establish a complete parallelism between the arithmetic properties of irreducible polynomials $F\in \kh[x]$,  encoded by their Okutsu frames, and the valuation-theoretic properties of their induced valuations $v_F$ on $\kh[x]$, encoded by their MLV chains.
In \cite{Rig}, it is shown that the natural restriction mapping $\ttt_{v^h}\to\ttt_v$ is an isomorphism of posets.

\section{Key polynomials over valued fields}\label{secKP}

In this section we introduce notation and well-known facts on key polynomials. Proofs and a more detailed exposition can be found in the survey \cite{KP}. 

For any field $L$, let $\irr(L)$ be the set of monic irreducible polynomials in $L[x]$. 

Let $(K,v)$ be a valued field. Let $k$ be the residue class field, $\g=v(K^*)$ the value group and $\gq=\g\otimes\Q$ the divisible hull of $\g$.

Let $\g\hk\Lambda$ be an extension of ordered abelian groups. We write simply $\Lambda\infty$ instead of $\Lambda\cup\{\infty\}$.
Consider a valuation on the polynomial ring $\kx$
$$
\mu\colon \kx\lra \Lambda\infty
$$
whose restriction to $K$ is $v$. 
Let $\p=\mu^{-1}(\infty)$ be the  support of $\mu$. 

The valuation $\mu$ induces in a natural way a valuation $\bar{\mu}$ on the field of fractions of $\kx/\p$; that is, $K(x)$ if $\p=0$, or $\kx/(f)$ if $\p=f\kx$ for some $f\in\irr(K)$.

The residue field $\km$ of $\mu$ is, by definition, the residue field of $\bar{\mu}$.

We say that $\mu$ is \emph{commensurable} (over $v$) if $\g_\mu/\g$ is a torsion group. In this case, there is a canonical embedding $\g_\mu\hookrightarrow \gq$. 
All valuations with nontrivial support are commensurable.

%\subsection{Basic properties of key polynomials}\label{subsecBPKP}
For any $\alpha\in\g_\mu$, consider the abelian groups:
$$
\ppa=\{g\in \kx\mid \mu(g)\ge \alpha\}\supset
\ppa^+=\{g\in \kx\mid \mu(g)> \alpha\}.
$$    
The \emph{graded algebra of $\mu$} is the integral domain:
$$
\ggm:=\gr_{\mu}(\kx)=\bigoplus\nolimits_{\alpha\in\g_\mu}\ppa/\ppa^+.
$$

There is a natural \emph{initial term} mapping $\inm\colon \kx\to \ggm$, given by $\inm\p=0$ and 
$$
\inm g= g+\pset_{\mu(g)}^+\in\pset_{\mu(g)}/\pset_{\mu(g)}^+, \qquad\mbox{if }\ g\in \kx\setminus\p.
$$

There is a natural embedding of graded algebras \ $\gg_v:=\op{gr}_v(K)\hookrightarrow \ggm$.

 The following definitions translate properties of the action of  $\mu$ on $\kx$ into algebraic relationships in the graded algebra $\ggm$.\e

\defn Let $g,\,h\in \kx$.

We say that $g,h$ are \emph{$\mu$-equivalent}, and we write $g\smu h$, if $\inm g=\inm h$. 

We say that $g$ is \emph{$\mu$-divisible} by $h$, and we write $h\mmu g$, if $\inm h\mid \inm g$ in $\ggm$.

We say that $g$ is $\mu$-irreducible if $(\inm g)\ggm$ is a nonzero prime ideal. 

We say that $g$ is $\mu$-minimal if $g\nmid_\mu f$ for all nonzero $f\in \kx$ with $\deg(f)<\deg(g)$.\e

The $\mu$-minimality condition admits a relevant characterization \cite[Prop. 2.3]{KP}.

\begin{lemma}\label{minimal0}
A polynomial  $g\in \kx\setminus K$ is $\mu$-minimal if and only if $\mu$ acts as follow on $g$-expansions:
$$f=\sum\nolimits_{0\le s}a_s g^s,\quad \deg(a_s)<\deg(g)\ \imp\ \mu(f)=\min\left\{\mu\left(a_sg^s\right)\mid 0\le s\right\}.$$
\end{lemma}

\defn A  \emph{(Maclane-Vaqui\'e) key polynomial} for $\mu$ is a monic polynomial in $\kx$ which is simultaneously  $\mu$-minimal and $\mu$-irreducible. 

The set of key polynomials for $\mu$ is denoted $\kpm$. \e

All $\phi\in\kpm$ are irreducible in $\kx$. Let $[\phi]_\mu\subset\kpm$ be the subset of all key polynomials  $\mu$-equivalent to $\phi$. Two $\mu$-equivalent key polynomials have the same degree \cite[Prop. 6.6]{KP}; hence, it makes sense to consider the degree $\deg\, [\phi]_\mu$ of a class.

The existence of key polynomials can be characterized as follows  \cite[Thm. 4.4]{KP}.

\begin{theorem}\label{KPempty}
%Let $\mu$ be a valuation on $\kx$, whose restriction to $K$ is $v$. 
The following conditions are equivalent.
\begin{enumerate}
\item $\kpm=\emptyset$.
\item $\mu$ is commensurable and $\km/k$ is an algebraic extension of fields.
\item $\ggm$ is a simple algebra (all nonzero homogeneous elements are units).
%\item $\ggm/\ggv$ is an algebraic extension of graded algebras.
\end{enumerate}
\end{theorem}

\defn
Suppose that  $\kpm\ne\emptyset$ and take $\phi\in\kpm$ of minimal degree. The \emph{degree} and \emph{singular value}  of $\mu$ are defined as
$$
\deg(\mu)=\deg(\phi),\qquad \sval(\mu)=\mu(\phi).
$$
The singular value is well defined by \cite[Thm. 3.9]{KP}. \e

%The \emph{weight} of $\mu$ is the value $$\wt(\mu)=\mu(f)/\deg(f)\in\La_\Q,$$ for all $\mu$-minimal $f\in\kx$. This value is independent of $f$ \cite[Thm. 3.9]{KP}. 

Another relevant invariant of a valuation $\mu$ on $\kx$ is its field $\ka=\kam$ of \emph{algebraic residues}, defined as the relative algebraic closure of $k$ in the residue field $\km$ of $\mu$.

Let $\d=\dm\subset\ggm$ be the subring of homogeneous elements of degree zero in the graded algebra. There are natural embeddings 
$$k\hk\ka\hk\d\hk\km.$$

The structure of $\d$ as a $\ka$-algebra plays an essential role in the description of the branches of a node in the valuative tree. %The following result is classical. Proofs can be found, for instance, in  \cite{KP}. 

\begin{theorem}\label{Delta}
Let $\mu$ be a valuation on $\kx$, whose restriction to $K$ is $v$. %Let $\mu\in\ttt$.
\begin{enumerate}
\item If $\kpm=\emptyset$, then $\ka=\d=\km$ is a countably generated extension of $k$.
\item If $\mu$ is incommensurable,  then $\ka=\d=\km$ is a finite extension of $k$.
\item If $\mu$ is commensurable and $\kpm\ne\emptyset$,  then $\d=\ka[\xi]$ and $\km=\ka(\xi)$, for some $\xi\in\d$ which is transcendental over $\ka$. 
\end{enumerate}
\end{theorem}

The valuations $\mu$ falling in case (3) of Theorem \ref{Delta} are said to be \emph{residually transcendental}.
There is a tight link between $\d$ and the set $\kpm$ \cite[Thm. 6.7]{KP}

\begin{theorem}\label{DeltaKP}
If $\kpm\ne\emptyset$, the residual ideal mapping $$\kpm\lra\op{Max}(\d),\qquad \phi\longmapsto \left(\inm(\phi)\ggm\right)\cap \d$$
induces a bijection between $\kpm/\!\smu$ and the maximal spectrum of $\d$.
\end{theorem}

If $\mu$ is incommensurable, then $\d$ is a field and $\op{Max}(\d)$ is a one-element set. In this case, $\kpm=[\phi]_\mu$, for any monic polynomial  $\phi\in\kx$ of minimal degree such that $\mu(\phi)$ is torsion-free over $\g$.

If $\mu$ is residually transcendental, Theorems \ref{Delta} and \ref{DeltaKP} yield a bijection %between 
$$\kpm/\!\smu\quad \longleftrightarrow\quad\op{Max}(\d)\quad \longleftrightarrow\quad \irr(\ka),$$ which depends on the choice of a generator $\xi$ for $\d$, as shown in Theorem \ref{Delta}(3).

\section{Tree of valuations with values in a fixed group}\label{secTreeLa}

Let $\ttt=\tla$ be the set of all valuations $\mu\colon \kx\to\La\infty$, whose restriction to $K$ is $v$.
This set admits a partial ordering. 
For $\mu,\nu\in \ttt$ we say that $\mu\le\nu$ if
$$
\mu(f)\le \nu(f),\qquad \forall\,f\in\kx.
$$
As usual, we write $\mu<\nu$ to indicate that $\mu\le\nu$ and $\mu\ne\nu$.

If $\mu\le\nu$, there is a canonical homomorphism of graded $\gg_v$-algebras:
$$\ggm\lra\ggn,\qquad \inm f\longmapsto
\begin{cases}\inu f,& \mbox{ if }\mu(f)=\nu(f),\\ 0,& \mbox{  if }\mu(f)<\nu(f).\end{cases}
$$

This poset $\ttt$ has the structure of a \emph{tree}. By this, we simply mean that all intervals 
$$
(-\infty,\mu\,]:=\left\{\rho\in\ttt\mid \rho\le\mu\right\}
$$
are totally ordered \cite[Thm. 2.4]{MLV}. \e

\defn A node $\mu\in\ttt$ is a \emph{leaf} if it  is a maximal element with respect to the ordering $\le$. Otherwise, we say that $\mu$ is an \emph{inner node}.

%It is easy to characterize when a node $\mu\in\ttt$ is a leaf in terms of algebraic properties of the valuation $\mu$ .

\begin{theorem}\cite[Thm. 2.3]{MLV}\label{maximal}
A node $\mu\in\ttt$ is a leaf if and only it $\kpm=\emptyset$.
\end{theorem}

All valuations with nontrivial support are leaves of $\ttt$, because they satisfy condition (2) of Theorem \ref{KPempty}. We call them \emph{finite leaves}.
The leaves of $\ttt$ having trivial support are \emph{valuation-algebraic} in Kuhlmann's terminology \cite{Kuhl}. We call them \emph{infinite leaves}.
We denote the set of leaves and subsets of finite and infinite leaves as follows:
$$
\ll(\ttt)=\lfin(\ttt)\sqcup\li(\ttt),
$$

\defn For a leaf $\mu\in\ll(\ttt)$ we define its \emph{degree} as:
$$
\deg(\mu)=\sup\left\{\deg(\rho)\mid \rho\in\ttt,\ \rho<\mu\right\}\in\N\infty.
$$\vskip0.2cm

A finite leaf  $\mu\in\lfin(\ttt)$  has $\supp(\mu)=f\kx$ for some monic irreducible  $f\in\kx$ and $\deg(\mu)=\deg(f)$. The infinite leaves may have finite or infinite degree.%\e

\subsection{Tangent directions and augmentations}\label{subsecTanDir}
Let $\mu,\,\nu$ be two nodes in $\ttt$ such that $\mu<\nu$.
Let $\tmn$ be the (nonempty) set of monic polynomials $\phi\in\kx$ of minimal degree satisfying $\mu(\phi)<\nu(\phi)$.

We say that $\tmn$ is the \emph{tangent direction} of $\mu$, determined by $\nu$. This terminology will be justified in section \ref{subsecTanSpace}, when we study the tangent space of $\ttt$.

The following properties of $\tmn$ were proven by Maclane for discrete rank-one valued fields, and generalized by Vaqui\'e to arbitrary valued fields \cite[Thm. 1.15]{Vaq}, \cite[Prop. 2.2, Cor. 2.5]{MLV}.

\begin{lemma}\label{propertiesTMN}
Let $\mu<\nu$ be two nodes in $\ttt$ and let $\phi\in\tmn$.
\begin{enumerate}
\item The polynomial $\phi$ belongs to $\kpm$ and $\tmn=[\phi]_\mu$. Also, $\deg(\mu)\le \deg(\nu)$.
\item For all nonzero $f\in\kx$ the equality $\mu(f)=\nu(f)$ holds if and only if $\phi\nmid_\mu f$.
%\item The kernel of the homomorphism $\ggm\to\ggn$ is the prime ideal $(\inm\phi)\ggm$.
\item If $\mu<\nu<\rho$ in $\ttt$, then $\tmr=\tmn$. In particular,
$$
\mu(f)=\rho(f)\sii \mu(f)=\nu(f),\qquad  \forall\,f\in\kx.
$$
%\item All nonzero elements in the image of $\,\ggm\to\ggn$ are units. 
\end{enumerate}
\end{lemma}

On the other hand, for any inner node $\mu\in\ttt$, all $\mu$-equivalence classes in $\kpm$ are the tangent direction of $\mu$ with respect to some $\nu\in\ttt$ such that $\mu<\nu$.  
Indeed, for any $\phi\in\kpm$ and any $\ga\in\Lambda\infty$ such that $\mu(\phi)<\ga$, we may construct the \emph{augmented valuation}  $\nu=[\mu;\phi,\ga]$,
defined in terms of $\phi$-expansions as
$$
f=\sum\nolimits_{0\le s}a_s\phi^s,\quad \deg(a_s)<\deg(\phi)\  \Longrightarrow\ \nu(f)=\min\{\mu(a_s)+s\ga\mid 0\le s\}.
$$
Note that $\nu(\phi)=\ga$.
The following properties of this augmented valuation  are also due to Maclane and Vaqui\'e \cite[Prop. 2.1]{MLV}.

\begin{lemma}\label{propertiesAug}
Let $\nu=[\mu;\phi,\ga]$ be an augmented valuation of $\mu$.
\begin{enumerate}
\item We have $\mu<\nu$ and $\,\tmn=[\phi]_\mu$. 
\item The value group of $\nu$ is $\gn=\gen{\g_{\mu,\phi},\ga}$, where $\g_{\mu,\phi}$ is the subgroup $$\g_{\mu,\phi}=\left\{\mu(a)\mid a\in\kx,\ 0\le \deg(a)<\deg(\phi)\right\}.$$
\item If $\ga=\infty$, then $\supp(\nu)=\phi\kx$.  If $\ga<\infty$, then $\phi$ is a key polynomial for $\nu$ of minimal degree. In both cases, $\deg(\nu)=\deg(\phi)$.
\end{enumerate}
\end{lemma}

For all $\phi_*\in\kpm$, $\ga_*\in\La\infty$ such that $\mu(\phi_*)<\ga_*$,  \cite[Lem. 2.8]{MLV} shows that
\begin{equation}\label{eqAug}
[\mu;\,\phi,\ga]=[\mu;\,\phi_*,\ga_*] \sii  \mu(\phi_*-\phi)\ge \ga=\ga_*\imp \phi\smu\phi_*.
\end{equation}

\subsection{The tangent space of $\ttt$}\label{subsecTanSpace}

For any inner node $\mu\in\ttt$, consider the quotient set
$$
\tdm=\left\{\nu\in\ttt\mid \mu<\nu\right\}/\!\tan,
$$
with respect to the equivalence relation $\tan$ which considers $\nu\tan\nu'$ if and only if $ (\mu,\nu]\cap(\mu,\nu']\ne\emptyset$.

The transitivity of $\tan$ follows easily from the fact that $\ttt$ is a tree. We denote by $[\nu]\stan$ the class of $\nu$. 
The elements of $\tdm$ can be identified with the tangent directions of $\mu$ defined in the last section.

\begin{proposition}\label{td=td}
For all inner nodes $\mu\in\ttt$, the association
$$
\phi\longmapsto t_\mu(\phi):=\left[\,[\mu;\,\phi,\ga]\,\right]\stan,\quad \ga\in\La\infty, \ \ga>\mu(\phi),
$$
is independent of the choice of $\ga$ and so it defines a mapping $t_\mu\colon\kpm\to\tdm$, which induces a bijection between  $\kpm/\!\smu$ and $\tdm$.
\end{proposition}

\begin{proof}
%By the very definition of $\tan$, the mapping $t_\mu$ is well defined. 
If $\mu(\phi)<\ga<\ga_*$, then $\nu<[\mu;\,\phi,\ga]<[\mu;\,\phi,\ga_*]$. Thus, $[\mu;\,\phi,\ga]\tan[\mu;\,\phi,\ga_*]$, so that the mapping $t_\mu$ is well defined.

Take $\nu \in \ttt$ such that $\mu<\nu$. For all $\phi\in\tmn$ we have
$$
\mu<[\mu;\,\phi,\nu(\phi)]\le\nu,
$$
by comparing their actions on $\phi$-expansions. Thus, $t_\mu(\phi)=[\nu]\stan$. This proves that $t_\mu$ is onto. 
Finally, let us show that, for all $\phi,\phi_*\in\kpm$, the equality $t_\mu(\phi)=t_\mu(\phi_*)$ holds if and only if $\phi\smu\phi_*$.

If $\phi\smu\phi_*$ and $\ga=\mu(\phi-\phi_*)>\mu(\phi)$, then $[\mu;\,\phi,\ga]=[\mu;\,\phi_*,\ga]$, by (\ref{eqAug}); thus $t_\mu(\phi)=t_\mu(\phi_*)$.
Conversely, if $t_\mu(\phi)=t_\mu(\phi_*)$, there exists $\rho\in\ttt$ such that
$$
\mu<\rho\le [\mu;\,\phi,\ga]\qquad\mbox{and}\qquad \mu<\rho\le [\mu;\,\phi_*,\ga_*],
$$
for some $\ga>\mu(\phi)$, $\ga_*>\mu(\phi_*)$. By \cite[Lem. 2.7]{MLV}, there exist $\dta,\dta_*\in\La$ such that
$[\mu;\,\phi,\dta]=\rho=[\mu;\,\phi_*,\dta_*]$. By (\ref{eqAug}), we have $\phi\smu\phi_*$. 
\end{proof}\e

%From now on, we identify both concepts. An element in any of the sets  $\kpm/\!\smu$, or $\tdm$, will be  simply called a \emph{tangent direction} of $\mu$.

By the remarks following Theorem \ref{DeltaKP}, $\tdm$ is a one-element set if $\mu$ is incommensurable, while there is a (non-canonical) bijection between $\tdm$ and $\irr(\kam)$, if $\mu$ is commensurable.\e

\defn The \emph{tangent space} of $\ttt$ is the set $\T(\ttt)$  containing all pairs $(\mu,t)$, where $\mu$ is an inner node in $\ttt$ and
$t\in\tdm$ is a tangent direction of $\mu$.

\section{Finite leaves}\label{secFinLeaves}
For any field $L$ and a monic irreducible polynomial $F\in\irr(L)$, we denote by $L_F$ the simple field extension $L[x]/(F)$. 

In this section, we assume that $\La$ contains the divisible closure of $\g$; that is, $\gq\subset\La$.
Under this assumption, the set $\lfin(\ttt)$ of finite leaves of $\ttt$ may be parametrized as
$$
\lfin(\ttt)=\left\{(\phi,\nb)\mid \phi\in\irr(K),\ \nb\ \mbox{ valuation on $K_\phi$ extending }v \right\},
$$
where we identify each pair $(\phi,\nb)$ with the following valuation with support $\phi\kx$:
$$
\nu\colon \kx\longtwoheadrightarrow K_\phi\stackrel{\nb}\lra\gq\infty
$$

%\nn{\bf Remark. }

Every simple field extension $L/K$ admits a finite number of extensions of $v$ to $L$.  
Any such extension determines an infinite number (if $K$ is infinite) of finite  leaves of $\ttt$, one for each  $\phi\in\irr(K)$ such that $K_\phi$ is $K$-isomorphic to $L$.
For instance, the valuation $v$ on $K$ determines the finite leaves $(x-a,v)$, for $a$ running in $K$.

%\subsection{Extensions of $v$ to simple field extensions of $K$}
Let us recall the description of all extensions of $v$ to simple finite extensions of $K$, which can be found (for instance) in \cite[Sec. 17]{endler}. 
Let us first describe all extensions of $v$ to an arbitrary algebraic extension $L$ of $K$. These extensions are commensurable over $v$; thus, we aim to describe the set:
$$
\ee(L)=\left\{w\colon L\lra \gq\infty\,\mid\, w\, \mbox{ valuation extending }v\right\}.
$$

Consider $K\subset \ks\subset\kb$, the separable closure of $K$ in a fixed algebraic closure $\kb$.

Let $\vb$ be a fixed extension of $v$ to $\kb$. %The group of values of $\vb$ is $\gq$ and its residue field is an algebraic closure of $k$. 
Let $K\subset \kh\subset \ks$ be the henselization of $K$ determined by the choice of $\vb$. Thus, $\kh$ is the fixed field of the decomposition group of the restriction of $\vb$ to $\ks$.

On the set  $\mon$ of all $K$-morphisms from $L$ to $\kb$, we define the following equivalence relation
$$
\la\sim_{\kh}\la' \ \sii\ \la'=\sg\circ\la\quad\mbox{for some}\quad \sg\in\aut(\kb/\kh).
$$
%We denote by $\mon{\kh}$ the quotient set $\mon{}/\!\sim_{\kh}$.

\begin{theorem}\label{endler}
	The  mapping $\mon\to\ee(L)$, defined by $\la\mapsto \vb\circ\la$, 
	induces a bijection between the quotient set $\mon{}/\!\sim_{\kh}$ and $\ee(L)$.
\end{theorem}

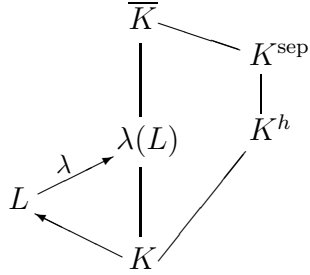
\begin{center}
\setlength{\unitlength}{4mm}
\begin{picture}(8,9)
\put(4,0){$K$}\put(3.6,4){$\la(L)$}\put(4,8){$\kb$}
\put(0,2){$L$}\put(8,4.3){$\kh$}\put(8,6.8){$\ks$}
\put(4.4,1.1){\line(0,1){2.3}}\put(4.4,5.1){\line(0,1){2.4}}\put(8.4,5.2){\line(0,1){1.3}}
\put(3.8,0.4){\vector(-2,1){2.9}}\put(1,2.5){\vector(2,1){2.5}}
\put(5,.3){\line(4,5){2.9}}\put(5,8.2){\line(3,-1){2.7}}
\put(1.6,3.2){\footnotesize{$\la$}}
\end{picture}
\end{center}\e

For instance, $\ee(\kb)$ is in bijection with  $\aut(\kb/\kh)\backslash\aut(\kb/K)$. Every right coset  $\aut(\kb/\kh)\,\sg$ determines the valuation $\vb\circ\sg$.

Suppose now that $L/K$ is a simple finite extension; that is, $L=K_\phi$ for some $\phi\in\irr(K)$. Since $\kh/K$ is a separable extension, the factorization of $\phi$ into a product of monic irreducible polynomials in $\kh[x]$ takes the form
$$
\phi=F_1\cdots F_r,
$$
with pairwise different $F_1,\dots,F_r\in\irr(\kh)$.
 Let $Z(\phi)\subset \kb$ be the set of zeros of $\phi$, avoiding multiplicities. We have a natural bijection
$$
Z(\phi)\lra\mon,\qquad \t\longmapsto \la_\t, 
$$
where $\la_\t$ is determined by $\la_\t\left(x+\phi\kx\right)=\t$.
Clearly,
$$
\begin{array}{rccl}
\la_\t \sim_{\kh} \la_{\t'}&\sii&\la_{\t'}=\sg\circ \la_\t&\quad\mbox{for some}\quad\sg\in\aut(\kb/\kh)\\
&\sii&\t'=\sg(\t)&\quad\mbox{for some}\quad\sg\in\aut(\kb/\kh)
\end{array}
$$
Therefore, $\la_\t \sim_{\kh} \la_{\t'}$ if and only if $\t$ and $\t'$ are roots of the same irreducible factor of $\phi$ over $\kh[x]$.

Let us choose an arbitrary root $\t_i\in Z(F_i)$ for each irreducible factor of $\phi$. By Theorem \ref{endler}, the set of valuations $\wb_{F_i}=\vb\circ\la_{\t_i}$ does not depend on the chosen roots and contains all extensions of $v$ to $L$.

\begin{theorem}\label{CorEndler}
There are $r$ extensions  of $v$ to $L=K_\phi$, given by $\wb_{F_1},\dots,\wb_{F_r}$. 
\end{theorem}

%Equivalently, if we denote $L_i=\kh_{F_i}=\kh[x]/F_i$ for $1\le i\le r$, we may define embeddings $L\hk L_i$ by using  the projections:$$L=\kx/\phi\hra \kh[x]/\phi\iso L_1\oplus\cdots \oplus L_r \lra  L_i.$$By the henselian property, the restriction of $\vb$ to $\kh$ has a unique extension $w_i$ to every $L_i/\kh$. It is easy to check that the restriction to $L$ of  $w_i$ coincides with the valuation $w_{F_i}$ for all $1\le i\le r$. 

%\subsection{Another parametrization of $\lfin(\ttt)$}
This description of the extensions of $v$ to simple finite extensions of $K$ yields a parametrization of the finite leaves by the set $\irr(\kh)$.
For all $F\in\irr(\kh)$ consider the finite leaf $w_F\in\lfin$ given by
$$
w_F(g)=\vb(g(\t))\quad \mbox{for all }g\in\kx,
$$
where $\t\in\kb$ is any root of $F$ in $\kb$. By the henselian property, this valuation $w_F$ is independent on the choice of $\t$.

Clearly, $\supp(w_F)=N(F)\kx$, where the ``norm" polynomial $N(F)\in\irr(K)$ is uniquely  determined by the equality $\left(F \khx\right)\cap \kx=N(F)\kx$.

Since the valuation induced by $w_F$ on $K_{N(F)}$ is precisely $\wb_F=\vb\circ\la_\t$, Theorem \ref{CorEndler} implies the following result.

\begin{theorem}\label{Lfin=Hensel}
If $\,\gq\subset \La$, we have a bijection 
$$
\irr(\kh)\lra \lfin(\ttt),\qquad F\longmapsto w_F=(N(F),\wb_F). 
$$
\end{theorem}

The inverse mapping associates to each pair $(\phi,\nb)\in\lfin(\ttt)$, the irreducible factor of $\phi$ over $\kh[x]$ canonically associated to $\nb$ by Theorem  \ref{CorEndler}.

%Note that  $w_F=(N(F),\wb_F)$, with respect to the initial parametrization of $\lfin$. 

\section{Infinite leaves and limit nodes}

In this section, we study the nodes of $\ttt$ which cannot be obtained by a finite chain of ordinary augmentations starting with  a degree-one valuation.
These nodes will be a kind of limit of certain totally ordered families of valuations in $\ttt$. 

\subsection{Totally ordered families of valuations}\label{subsecTOF}
Consider a totally ordered family of inner nodes of $\ttt$, not containing a maximal element:
$$
\aa=\cfa,\quad \rho_i\in\ttt.
$$
We assume that $\aa$ is parameterized by a totally ordered set set $A$ of indices such that the mapping $A\to\aa$ determined by $i\mapsto \ri$ is an isomorphism of totally ordered sets.

By Lemma \ref{propertiesTMN}, the degree function $\deg\colon \aa\to\N$
is order-preserving. Hence, these families fall into two radically different cases:\e
 
 (a) \ The set $\deg(\aa)$ is unbounded in $\N$. We say that $\aa$ has \emph{unbounded degree}.\e
 
 (b) \ There exists $i_0\in A$ such that $\deg(\ri)=\deg(\rho_{i_0})$ for all $i\ge i_0$. 
We say that $\aa$ is a \emph{continuous family} of stable degree $m(\aa)=\deg(\rho_{i_0})$.\e

%\subsubsection*{Stability of polynomials}\mbox{\null}

In any case, $\aa$ determines a unique tangent direction of every valuation $\ri\in\aa$. Indeed, Lemma \ref{propertiesTMN} shows that $\ty(\ri,\rj)=\ty(\ri,\rk)$ for all $i<j<k$ in $A$.

We denote by $\ty(\ri,\aa)$ this common tangent direction. By Lemma \ref{propertiesTMN}, there exists a key polynomial $\varphi_i\in\kp(\ri)$ such that $\ty(\ri,\aa)=[\varphi_i]_{\ri}$,
and for any nonzero polynomial $f\in\kx$ we have
$$
\varphi_i\mid_{\ri}f \ \sii\ \ri(f)<\rj(f)\ \mbox{ for all }\ j>i \ \mbox{ in }A.
$$
$$
\varphi_i\nmid_{\ri}f \ \sii\ \ri(f)=\rj(f)\ \mbox{ for all }\ j>i \ \mbox{ in }A.
$$\vskip0.2cm

\defn 
We say that a nonzero $f\in\kx$ is \emph{$\aa$-stable} if it satisfies $\varphi_i\nmid_{\ri}f$, for some index $i\in A$.
In this case, we denote its stable value by  $$\rha(f)=\max\{\ri(f)\mid i\in A\}.$$ 
We obtain in this way a \emph{stability function} 
$\rha\colon \kx_\aa\to \La\infty$, defined only on the subset  $\kx_\aa\subset \kx$ formed by the $\aa$-stable polynomials.\e

%Note that $f$ is \emph{$\aa$-unstable} if and only if $\ri(f)<\rj(f)$  for all pairs $i<j$ in $A$; or, equivalently, $\varphi_i\mid_{\ri}f$ for all $i\in A$. 

The following basic properties of the function $\rha$ are obvious:
\begin{itemize}
\item $(\rha)_{\mid K}=v$,
\item $f,g\in \kx_\aa\ \imp\ fg\in\kx_\aa\ \mbox{ and } \ \rha(fg)=\rha(f)+\rha(g)$,
\item $f,g,f+g\in \kx_\aa\ \imp\ \rha(f+g)\ge\min\{\rha(f),\rha(g)\}$.
\end{itemize}

In particular, if $\kx_\aa=\kx$, then the  function $\rha$ is a valuation in $\ttt$.\e

\defn
If all the polynomials in $\kx$ are $\aa$-stable, we say that the valuation $\rha$ is the \emph{stable limit} of $\aa$. 
In this case, we write $\rha=\lim(\aa)=\lim_{i\in A}\ri$.\e

\begin{proposition}\cite[Prop. 3.1]{MLV}\label{leaves}
If $\aa=\cfa$ has a stable limit, then $\rha$ has trivial support and $\kp(\rha)=\emptyset$. In particular, $\rha$ is an infinite leaf of the tree $\ttt$.
\end{proposition}

Let us see a necessary condition for a polynomial to be $\aa$-unstable. 

\begin{lemma}\label{aunstab}
All $\aa$-unstable polynomials $f$ satisfy $\deg(f)\ge \deg(\ri)$ for all $i\in A$.
\end{lemma}

\begin{proof}
Let $\ty(\ri,\aa)=[\varphi_i]_{\ri}$ for some $i\in A$. If $f\in\kx$ has  $\deg(f)<\deg(\ri)$, then $\deg(f)<\deg(\varphi_i)$. Hence, $\varphi_i\nmid_{\ri}f$, contradicting the unstability of $f$.   
\end{proof}

\begin{corollary}\label{unbIsStab}
Every totally ordered family of unbounded degree has a stable limit. %, and the function $$\sval\colon \aa\lra\La,\qquad \ri\longmapsto \sval(\ri) $$is unbounded too.
\end{corollary}

%The following alternative characterization of $\aa$-stability is very useful in practice. It follows immediately from  \cite[Cor. 2.5,(2)]{MLV}.

%\begin{lemma}\label{astab} Let $\aa=\cfa$ be a totally ordered family and let $f\in \kx$. Then,$$\ri(f)=\rj(f)  \mbox{ for some $i<j$ in }A\ \imp \  \ri(f)=\rk(f)  \mbox{ for all $i\le k$ in }A.$$In particular, $\rha(f)=\ri(f)$.\end{lemma}

\subsection{Continuous families and limit augmentations}\label{subsecCont}

Let $\cc=\cfa$ be a \emph{continuous family} of valuations in $\ttt$  of stable degree $m=m(\cc)$.

%Let $\aa=\cfa$ be a totally ordered family of valuations in $\ttt$ (which does not contain a maximal element) of stable degree $m$. 

The following result shows that the valuations of maximal degree $m$ in $\cc$ are ``close" one to another in a  certain sense. 

\begin{lemma}\label{samedegree}
Let $\mu<\nu$ be two inner nodes of $\ttt$ of degree $m$. Take any $\phi\in\kpn$ of degree $m$. Then, $\phi\in\kpm$ and  $\nu= [\mu;\,\phi, \nu(\phi)]$.    
\end{lemma}

\begin{proof}
For all $a \in\kx$ with $\deg(a)<m$, we have $\mu(a)=\nu(a)$, because any $\varphi\in \ty(\mu,\nu)$ satisfies $\varphi\nmid_\mu a$.

Necessarily $\mu(\phi)<\nu(\phi)$, because the equality $\mu(\phi)=\nu(\phi)$ leads to a contradiction. Indeed, for all $f\in\kx$ with $\phi$-expansion $f=\sum_{0\le s}a_s\phi^s$, we would have
$$
\mu(f)\ge\min_{0\le s}\{\mu\left(a_s\phi^s\right)\}=\min_{0\le s}\{\nu\left(a_s\phi^s\right)\}=\nu(f),
$$
showing that $\mu\ge\nu$, against our assumption.

By Lemma \ref{propertiesTMN}, $\phi$ is a key polynomial for $\mu$ such that $\ty(\mu,\nu)=[\phi]_\mu$. The equality $\nu= [\mu;\,\phi, \nu(\phi)]$ follows then from the action of both valuations on $\phi$-expansions.
\end{proof}%\e

\subsubsection{Stable group of a continuous family}\label{subsubsecG0} \mbox{\null}

Let $\rho\in\ttt$ be an inner node. By Lemma \ref{minimal0}, applied to any key polynomial for $\rho$ of minimal degree,  we have $\grh=\gen{\grh^0,\sval(\rho)}$, where $\grh^0$ is the subgroup:
$$
\grh^0=\left\{\rho(a)\mid a\in\kx,\ 0\le\deg(a)<\deg(\rho)\right\}.
$$

Moreover, $\grh^0$ is commensurable over $\g$. Indeed, for all $a\in\kx$ of degree less than $\deg(\rho)$, the initial term $\inn_\rho a\in\gg_\rho$ is algebraic over the graded algebra $\gg_v$ \cite[Prop. 3.5]{KP}; thus, $n\rho(a)$ belongs to $\g$ for some $n\in\N$.

The index $\erel(\rho)=\left(\grh\colon \grh^0\right)$ is called the \emph{relative ramification index} of $\rho$. 

If $\rho$ is incommensurable over $v$, then $\erel(\rho)=\infty$.
%Clearly, the property $\erel(\rho)=1$ is equivalent to $\sval(\rho)\in\grh^0$.

Let us denote the set of stable values of nonzero $\cc$-stable polynomials by
$$\gc=\rhc\left(\kx_\cc\setminus\{0\}\right)\subset\La.$$ 
The following lemma shows that $\gc$ is a subgroup of $\La$. It is called 
the \emph{stable group} of the continuous family $\cc$. 

\begin{lemma}\label{GAcomm}
Let $\cc=\cfa$ be a \emph{continuous family} in $\ttt$ of stable degree $m$. Then, $\g_{\ri}^0=\gc$ \,for all $\ri\in \cc$ of degree $m$.
\end{lemma}

\begin{proof}
Suppose that $\deg(\ri)=m$ for some $i\in A$.
Since $\cc$ contains no maximal element, there exists $j\in A$ such that $i<j$. Since the degree function preserves the ordering, we have $\deg(\ri)=\deg(\rj)=m$.  
By Lemma \ref{samedegree}, $\ty(\ri,\rj)=[\varphi]_{\ri}$, for any $\varphi\in\kp(\rj)$ of minimal degree $m$.

Now, for all $a\in\kx$ with $0\le\deg(a)<m$, we have $\ri(a)=\rj(a)=\rhc(a)$ because $\varphi\nmid_{\ri}a$. This proves $\g_{\ri}^0=\g_{\rj}^0\subset \gc$ already.

Let us show the inclusion $\gc\subset \g_{\ri}^0$. Let $\ga=\rhc(f)\in\gc$ be the stable value of a nonzero $\cc$-stable  $f\in\kx$. Let $i\le k<\ell$ in $A$ such that $\deg(\rk)=\deg(\rl)=m$ and $\rk(f)=\rl(f)$. By \cite[Cor. 2.6]{MLV}, the element $\inn_{\rl} f$ is a unit. By \cite[Prop. 3.5]{KP}, there exists $a\in\kx$ of degree less than $m$ such that $\inn_{\rl} f=\inn_{\rl} a$. In particular, $\rl(f)=\rl(a)$ belongs to $\g_{\rl}^0=\g_{\ri}^0$.
\end{proof}%\e

\subsubsection{Limit key polynomials and limit augmentations}\label{subsubsecLimit}\mbox{\null}

The set $\kpi(\cc)$ of \emph{limit key polynomials} for $\cc$ is the set  of all monic $\cc$-unstable polynomials of minimal degree. Since the product of stable polynomials is stable, all limit key polynomials are irreducible in $ \kx$.

The \emph{unstable degree} $\mi=\mi(\cc)$ is the common degree of all limit key polynomials for $\cc$.
 If all polynomials in $\kx$ are $\cc$-stable, then $\kpi(\cc)=\emptyset$ and we agree that $\mi=\infty$. 
 By Lemma \ref{aunstab}, $\mi\ge m$.

Take any limit key polynomial $\phi\in\kpi\left(\cc\right)$, and choose $\ga\in\La\infty$ such that $$\ri(\phi)<\ga\quad \mbox{for all } \ i\in A.$$
The \emph{limit augmentation}   $\nu=[\cc;\phi,\ga]$ acts as follows on $\phi$-expansions: 
$$f=\sum\nolimits_{0\le s}a_s\phi^s,\quad\deg(a_s)<\mi\ \imp\ \nu(f)=\min\{\rhc(a_s)+s\ga\mid 0\le s\}.
$$
%Since $\deg(a_s)<\mi$, all coefficients $a_s$ are $\cc$-stable. %Note that $\nu(\phi)=\ga$.
The following properties of $\nu$ can be found in 
\cite[Sec. 1.4]{Vaq} and \cite[Sec. 7]{KP}.

\begin{proposition}\label{extensionlim} 
The augmentation $\nu=[\cc;\phi,\ga]$ is a valuation in $\ttt$  such that $\ri<\nu$ for all $i\in A$.
If $\ga<\infty$, then $\nu$ is an inner node, $\gn=\gen{\gc,\ga}$ and $\phi$ is a key polynomial for $\nu$, of minimal degree; thus, $\deg(\nu)=\mi$ and $\,\sval(\nu)=\ga$.  

If $\ga=\infty$, then $\gn=\gc$ and the support of $\nu$ is $\phi\kx$; thus, $\nu$ is a finite leaf of $\ttt$.   
\end{proposition}

If $\mi=m$, then for all $i\in A$ with maximal degree $\deg(\ri)=m$, Lemma \ref{samedegree} shows that
$\nu=[\ri;\,\phi,\nu(\phi)]$.
Thus, any limit augmentation of $\cc$ is equal to an ordinary augmentation of some $\ri$. 

Therefore, we may discard the families $\cc$ of stable degree $m=\mi$ because they do not contribute to produce new nodes in $\ttt$ by limit augmentations.  

%\begin{definition}\label{defContFam}\rm  \end{definition}

Summing up, a continuous family $\cc$ of valuations has three possibilities:
 
 (a) \ If $m<\mi=\infty$, then $\cc$ has a stable limit.
 
 (b) \ If $m=\mi<\infty$, we say that $\cc$ is \emph{inessential}.

 (c) \ If $m<\mi<\infty$, we say that $\cc$ is \emph{essential}.\e

The continuous families having a stable limit determine infinite leaves of the tree $\ttt$, by Proposition \ref{leaves}. 
The essential continuous families determine inner \emph{limit nodes} of $\ttt$ as limit augmentations of the family. 
By the uniqueness of \mlv chains \cite[Thm. 4.7]{MLV}, these limit nodes cannot be obtained by chains of ordinary augmentations starting with valuations that are smaller than some $\ri\in\cc$. 

The limit key polynomials are an essential ingredient in the construction of these limit nodes. Let us describe them in more detail.

\begin{lemma}\label{allLKP}
Let $\cc$ be an essential continuous family and let $\phi\in\kpi(\cc)$. Then,
$$
\kpi(\cc)=\left\{\phi+a\mid a\in\kx,\ \deg(a)<\mi, \ \rhc(a)>\rho_i(\phi) \ \mbox{ for all }i\in A\right\}.
$$
\end{lemma}

\begin{proof}
Let $\varphi\in\kx$ be a monic polynomial of degree $\mi$. Then,
$\varphi=\phi+a$ for some $a\in\kx$ with $\deg(a)<\mi$.
Since $a$ is $\cc$-stable, there exists $i_0\in A$ such that
$$
\rhc(a)=\rho_j(a)\quad\mbox{ for all }\  i_0\le j.
$$

For all indices $i_0\le j<k$ we have $\rho_j(\phi)<\rho_k(\phi)$ and $\rho_j(a)=\rho_k(a)=\rhc(a)$. 

Suppose that $\rhc(a)>\rho_i(\phi)$ for all $i\in A$. From $\rhc(a)>\rho_k(\phi)$ we deduce that $\rho_j(\varphi)=\rho_j(\phi)<\rho_k(\phi)=\rho_k(\varphi)$. Thus, $\varphi$ belongs to $\kpi(\cc)$.

Suppose that $\rhc(a)\le \ri(\phi)$ for some $i\in A$. Then, for all $\max\{i_0,i\}<j<k$ we have $\rho_j(\varphi)=\rhc(a)=\rho_k(\varphi)$. Thus, $\varphi$ is $\cc$-stable.
\end{proof}

\subsection{Equivalence of totally ordered families of valuations}\label{subsecEquivCont}

Let $\aa=\cfa$, $\bb=\cfb$ be totally ordered families in $\ttt$, not containing a maximal element. 

We say that $\aa$ and $\bb$ are \emph{equivalent}, and we write $\aa\sim\bb$, if they are cofinal in each other.

Obviously, two equivalent families either both have unbounded degree, or both have stable degree. Also, in the latter case they have the same stable and unstable degrees, and the same stable group.

Two totally ordered families in the same class have the same limit behavior.

\begin{proposition}\label{equiv=lim}
Let $\aa$, $\bb$ be totally ordered families in $\ttt$, not containing a maximal element. Suppose that none of them is a continuous inessential family.  Then,  
$$
\aa\sim\bb \ \sii \ \kx_\aa=\kx_\bb \ \mbox{ and }\ \rha=\rhb.
$$
\end{proposition}

\begin{proof}
Suppose $\aa\sim\bb$. Take any $f\in\kx_\aa$, so that there exists $i_0\in A$ such that
$$
\rha(f)=\ri(f)\quad \mbox{for all }\ i\ge i_0. 
$$
Since $\bb$ is cofinal in $\aa$,  there exists $j_0\in B$ such that $\rho_{i_0}<\zeta_{j_0}$. Take any $j\in B$, $j>j_0$; since  $\aa$ is cofinal in $\bb$,  there exists $i\in A$ such that
$$
\rho_{i_0}<\zeta_{j_0}<\zj<\ri.
$$
Necessarily, $\rho_{i_0}(f)=\zeta_{j_0}(f)=\zj(f)=\rhb(f)$. Hence, $f$ is $\bb$-stable and $\rha(f)=\rhb(f)$.

The symmetry of the argument shows that $\kx_\aa=\kx_\bb$ and $\rha=\rhb$.

Conversely, suppose that $\kx_\aa=\kx_\bb$ and $\rha=\rhb$. 
Let us first assume that both families have a stable limit; that is,  $\kx_\aa=\kx_\bb=\kx$.

Since all valuations $\ri, \zj$ are bounded above by the valuation $\rha=\rhb$, the set $\aa\cup\bb$ is totally ordered. If $\aa$ and $\bb$ were not cofinal in each other, there would exist (for instance) some $j\in B$ such that $\zj>\ri$ for all $i\in A$.  But this implies $\zj\ge\rha=\rhb$, leading to $\zk>\rhb$ for all $k>j$ in $B$. This is a contradiction.

Suppose now $\kx_\aa=\kx_\bb\subsetneq\kx$. Take any $\phi\in\kpi(\aa)=\kpi(\bb)$. Since $\rha=\rhb$, the following limit augmentations coincide
$$
\nu:=[\aa;\,\phi, \infty]=[\bb;\,\phi, \infty],
$$
because they have the same action on $\phi$-expansions. As above, the set $\aa\cup\bb$ is totally ordered  and $\aa$ and $\bb$ are
cofinal in each other, unless there exists (for instance) some $j\in B$ such that $\zj>\ri$ for all $i\in A$. This leads to a contradiction too.

Indeed, let $\ty(\zj,\zk)=[\varphi]_{\zj}$ for some $k\in B$ such that $k>j$. Since $\zj(\varphi)<\zk(\varphi)$, the augmentation $\mu=[\zj;\,\varphi,\zk(\varphi)]$ satisfies $\zj<\mu\le\zk$. %By taking $j,k$ large enough we may assume $\deg(\zj)=\deg(\zk)=m(\bb)$. 
Since  
$$
\deg(\varphi)=\deg(\mu)\le\deg(\zk)\le m(\bb)<\mi(\bb),
$$
we deduce that $\varphi$ is $\bb$-stable. By our hypothesis, $\varphi$ is $\aa$-stable too, and this implies $\zj(\varphi)=\rha(\varphi)=\zk(\varphi)$, which is a contradiction.
\end{proof}\e

\begin{corollary}\label{sameLKP}
Let $\cc$, $\cc'$ be two equivalent essential continuous families in $\ttt$. Then, $\kpi(\cc)=\kpi(\cc')$.  
\end{corollary}

Since all totally ordered sets admit cofinal well-ordered subsets, in every class of totally ordered families of valuations in $\ttt$ there are well-ordered families.

In this vein, for  any given class, we want to study the existence of ``nice" families in the class, having special properties.

%\begin{lemma}\label{numerable}Let $\aa=\cfa$ be a totally ordered family of unbounded degree. Then, it admits a countable cofinal subfamily $\bb=\left(\rho_{i_n}\right)_{n\in\N}$ such that all $\rho_i$ are commensurable over $v$ and  $\deg(\rho_{i_m})<\deg(\rho_{i_n})$ for all $m<n$. \end{lemma}

%\begin{proof}Take any $i_0\in A$. We construct $\bb$ by a recurrent procedure. For all $n\in\N$ we consider any choice $i_{n}\in A$ such that $\deg(\rho_{i_{n}})>\deg(\rho_{i_{n-1}})$.   Since the degree mapping preserves the order, it is obvious that this subfamiliy $\bb$ is cofinal in $\aa$.   \end{proof}

\begin{lemma}\label{numerable}
Let $\aa=\cfa$ be a totally ordered family of unbounded degree. Then, there is a countable equivalent family $\bb=\left(\zeta_n\right)_{n\in \N}$ such that all $\zeta_n$ are commensurable over $v$ and  $\deg(\zeta_{m})<\deg(\zeta_{n})$ for all $m<n$. 
\end{lemma}

\begin{proof}
By Corollary \ref{unbIsStab}, $\aa$ has a stable limit $\mu=\lim_{i\in A}\rho_i$. The Maclane--Vaqui\'e theorem (Theorem \ref{main}) shows that $\mu$ is the stable limit of a countable infinite Maclane--Vaqui\'e  chain   $\bb=\left(\zeta_n\right)_{n\in \N}$ such that all $\zeta_n$ are commensurable over $v$ and  $\deg(\zeta_{m})<\deg(\zeta_{n})$ for all $m<n$. 
By Proposition \ref{equiv=lim}, $\aa$ and $\bb$ are equivalent. 
\end{proof}

\begin{lemma}\label{specialCont}
Let $\cc=\cfa$ be a continuous family such that $m(\cc)<\mi(\cc)$. Then, there is an equivalent family $\bb=\left(\zj\right)_{j\in B}$ satisfying the following properties:
\begin{enumerate}
\item $\bb$ is well-ordered and all valuations $\zi\in\bb$ have degree $m(\cc)$.
\item All valuations $\zj\in\bb$ have relative ramification index equal to one.
\item All valuations $\zj\in\bb$ have value group $\g_{\zj}=\gb$.
\item For all $j\in B$, there exist $\chi_j\in\kp(\zj)$ of minimal degree, such that
$$
j<k\ \mbox{in }B \ \imp\ \chi_j\not\sim_{\zj}\chi_k\ \mbox{ and }\ \zk=[\zj;\,\chi_k,\sval(\zk)].
$$
\end{enumerate}
\end{lemma}

\begin{proof}
By replacing $\cc$ with a cofinal subfamily, we may assume that $\cc$ is well-ordered and  $\deg(\ri)=m(\cc)$ for all $i\in A$.  By Lemma \ref{GAcomm}, $\g_{\ri}^0=\gc$ for all $i\in A$.

Let us construct an equivalent family $\bb=\left(\zj\right)_{j\in B}$ all whose valuations have degree $m(\cc)$, are commensurable and have relative ramification index equal to one; that is, $\sval(\zj)\in \gc=\gb$ for all $j\in B$.

Let $\ccb$ be the subset of $\cc$ formed by all $\ri$ such that $\sval(\ri)\not\in \gc$. 
Let us first construct, for each $\ri\in\ccb$, a valuation $\rho'_i$ such that $\sval(\rho'_i)\in\gc$ and $\ri<\rho'_i\le \rk$ for some $k$ in $A$.

For any given $\ri\in\ccb$, let $\varphi\in\kp(\ri)$ be a key polynomial of minimal degree. Then, $\sval(\ri)=\ri(\varphi)\not\in\gc$.
Since $\deg(\varphi)=\deg(\ri)<\mi(\cc)$, the polynomial $\varphi$ is $\cc$-stable. Take $i<j<k$ in $A$, such that $\rj(\varphi)=\rk(\varphi)=\rhc(\varphi)$. Since $\rhc(\varphi)\in \gc$, the inequality $\ri(\varphi)< \rk(\varphi)$ must be strict.

Now, take $\chi\in\kp(\rk)$ a key polynomial of minimal degree. By \cite[Thm. 3.9]{KP}, we have $\rk(\varphi)\le\rk(\chi)=\sval(\rk)$. By Lemma \ref{samedegree},
$\rk=[\ri;\, \chi,\sval(\rk)]$. Since $\ri(\varphi)< \rk(\varphi)\le \rk(\chi)$, the augmented valuation
$$
\rho'_i:=[\ri;\,\chi,\rk(\varphi)]
$$
satisfies $\sval(\rho'_i)=\rho'_i(\chi)=\rk(\varphi)\in\gc$ and $\ri<\rho'_i\le \rk$.

Consider the totally ordered family $\cc'=\cc\cup \{\rho'_i\mid \ri\in\ccb\}$. Obviously,  $\cc$ and $\cc'$ are equivalent. Also, by construction, the subfamily $\bb$ of $\cc'$ formed by all valuations $\rho$ such that $\sval(\rho)\in\gc$ is cofinal in $\cc'$. Thus, the family $\bb$ is equivalent to $\cc$ and satisfies conditions (1) and (2). Since $\g_{\zj}=\gc=\g_{\zj}^0$ for all $\zj\in\bb$, condition (3) follows from Lemma \ref{GAcomm}.

Let us prove that condition (4) holds too, for this family $\bb$.
Let us choose arbitrary key polynomials $\chi_j\in\kp(\zj)$ of minimal degree, for all valuations $\zj\in\bb$. By Lemma \ref{samedegree},  $\zk=[\zj;\,\chi_k,\sval(\zk)]$ for all $j<k$ in $B$. 

The class $\ty(\zj,\bb)=\ty(\zj,\zk)=[\chi_k]_{\zj}$ depends only on $\zj$ and $\bb$. 
Thus, the condition $\chi_j\not\sim_{\zj}\chi_k$, which is equivalent to $\chi_j\not\in \ty(\zj,\bb)$, depends only on $\chi_j$.

Let us show that there exists a key polynomial $\chi'_j\in\kp(\zj)$ of minimal degree, satisfying $\chi'_j\not\in \ty(\zj,\bb)$. Since $\zj$ has relative ramification index equal to one, we have $\sval(\zj)\in\g_{\zj}^0$ and there exists $a\in\kx$ of degree less than $m$ such that $\zj(a)=\sval(\zj)=\zj(\chi_j)$.  By \cite[Prop. 6.3]{KP}, $\chi'_j=\chi_j+a$ is a key polynomial  for $\zj$ of minimal degree $m$, such that $\chi'_j\not\sim_{\zj}\chi_j$. 
Therefore, at least one of the two key polynomials $\chi_j$, $\chi'_j$ does not fall in the class $\ty(\zj,\bb)$.
\end{proof}

 \section{Paths in the tree $\ttt$}\label{secCtDepth}
 
 \subsection{\mlv chains}\label{subsecMLV}

In this section, we review the fundamental theorem of \mlv describing how to reach all nodes in $\ttt$ by a combination of ordinary augmentations, limit augmentations and stable limits \cite{mcla,Vaq}. All results are extracted from the survey \cite{MLV}.

Consider a finite, or countably infinite, chain of nodes in $\ttt$
\begin{equation}\label{depthMLV}
\mu_0\ \stackrel{\phi_1,\ga_1}\lra\  \mu_1\ \stackrel{\phi_2,\ga_2}\lra\ \cdots
\ \lra\ \mu_{n-1} 
\ \stackrel{\phi_{n},\ga_{n}}\lra\ \mu_{n} \ \lra\ \cdots
\end{equation}
in which every node is an augmentation  of the previous node, of one of the following two types:\e

\emph{Ordinary augmentation}: \ $\mu_{n+1}=[\mu_n;\, \phi_{n+1},\ga_{n+1}]$, for some $\phi_{n+1}\in\kp(\mu_n)$.\e

\emph{Limit augmentation}:  \ $\mu_{n+1}=[\aa_n;\, \phi_{n+1},\ga_{n+1}]$,  for some $\phi_{n+1}\in\kpi(\aa_n)$, where $\aa_n$ is an essential continuous family whose first valuation is $\mu_n$.\e

We consider an implicit choice of a key  polynomial $\phi_0\in\kp(\mu_0)$ of minimal degree, and we denote $\ga_0=\mu_0(\phi_0)$.

%All values $\ga_n$ belong to $\La\infty$ and satisfy $\ga_n=\sval(\mu_n)<\ga_{n+1}$. 

Therefore, for all $n$ such that $\mu_n$ is an inner node of $\ttt$, the polynomial $\phi_n$ is a key polynomial for $\mu_n$ of minimal degree, and we have
$$%\begin{equation}\label{mgan}
m_n:=\deg(\mu_n)=\deg(\phi_n),\qquad \sval(\mu_n)=\mu_n(\phi_n)=\ga_n.
$$%\end{equation}\vskip0.2cm

\defn
A chain of mixed augmentations as in (\ref{depthMLV}) is said to be  a \emph{\mlv (abbreviated MLV) chain}  if $\deg(\mu_0)=1$ and every augmentation step satisfies:
\begin{itemize}
\item If $\,\mu_n\to\mu_{n+1}\,$ is ordinary, then $\ \deg(\mu_n)<\deg\,\ty(\mu_n,\mu_{n+1})$.
\item If $\,\mu_n\to\mu_{n+1}\,$ is limit, then $\ \deg(\mu_n)=m(\aa_n)$ and $\ \phi_n\not \in\ty(\mu_n,\mu_{n+1})$. 
\end{itemize}\e

Let $0\le r\le \infty$ be the length of a MLV chain. For $n<r$, all nodes $\mu_n$ are residually transcendental valuations. Indeed, in all augmentations of the chain, either ordinary or limit, we have
$$
\phi_n\in\kp(\mu_n),\quad \phi_n\not\in\ty(\mu_n,\mu_{n+1}).
$$
Hence, $\kp(\mu_n)$ contains at least two different $\mu_n$-equivalence classes. By the remarks at the end of Section \ref{secKP}, this implies that $\mu_n$ is commensurable. 

\begin{theorem}\label{main}
Every node $\nu\in\ttt$ falls in one, and only one, of the following cases.  \e

(a) \ It is the last valuation of a finite MLV chain.
$$ \mu_0\ \stackrel{\phi_1,\ga_1}\lra\  \mu_1\ \stackrel{\phi_2,\ga_2}\lra\ \cdots\ \lra\ \mu_{r-1}\ \stackrel{\phi_{r},\ga_{r}}\lra\ \mu_{r}=\nu.$$

(b) \ It is the stable limit of a continuous family $\aa=\cfa$ of augmentations whose first valuation $\mu_r$ falls in case (a):
$$ \mu_0\ \stackrel{\phi_1,\ga_1}\lra\ \mu_1\ \stackrel{\phi_2,\ga_2}\lra\   \cdots\ \lra\ \mu_{r-1}\ \stackrel{\phi_{r},\ga_{r}}\lra\ \mu_{r}\  \stackrel{\cfa}\lra\  \rha=\nu.
$$
Moreover, we may assume that $\deg(\mu_r)=m(\aa)$ and  $\phi_r\not\in\ty(\mu_r,\nu)$.\e

(c) \ It is the stable limit, $\nu=\lim_{n\in\N}\,\mu_n$, of an  infinite MLV chain.
\end{theorem}

The inner nodes and the finite leaves of $\ttt$ fall in case (a). These are the ``bien-specifi\'ees" valuations in Vaqui\'e's terminology. 

We denote by  $\lci(\ttt),\,\lui(\ttt)\subset\li(\ttt)$ the subsets of infinite leaves falling in cases (b), (c), respectively. The infinite leaves in $\lci(\ttt)$ have finite degree and those in $\lui(\ttt)$ have infinite degree.

Also, Lemma \ref{propertiesTMN} shows that in all cases displayed in Theorem \ref{main}, we have
$$%\begin{equation}\label{ga}
\phi_n\not\in\ty(\mu_n,\nu) \quad \mbox{ and }\quad \nu(\phi_n)=\mu_n(\phi_n)=\ga_n=\sval(\mu_n)\quad \mbox{for all }\ n.
$$%\end{equation}

The main advantage of MLV chains is that their nodes are essentially unique, so that we may read in them several data intrinsically associated to the valuation $\mu$.

For instance the sequence $(m_n)_{n\ge0}$ and the character ``ordinary" or ``limit" of each augmentation step $\mu_n\to\mu_{n+1}$ are intrisic features of $\nu$ \cite[Thm. 4.7]{MLV}.

Thus, we may define order preserving functions
$$
\dep,\ldp\colon\ttt\lra\N\infty,
$$
where $\dep(\nu)$ is the length of the MLV chain  underlying $\nu$, and
 $\ldp(\nu)$ counts the number of limit augmentations in this MLV chain
\footnote{Thus, all valuations of both types (a) and (b) have a finite depth. At this point, we are not following the convention of \cite{MLV}, where the valuations of type (a) were said to have \emph{finite depth} while those of type (b) were said to have \emph{quasi-finite depth}}.  

The arguments in the proof of \cite[Lem. 4.2]{MLV} show that these functions preserve the ordering.

\subsection{Decoding  MLV chains for arithmetic and geometric applications}
\label{subsecAMdata}\mbox{\null}

Besides their intrinsic theoretical interest, MLV chains encode a large amount of information which can
be useful in several contexts. In this section, we  describe  a concrete MLV chain in full
generality, and then we  interpret it from both the number theoretic and the geometric
perspective. In the former case, we will see how to describe the decomposition of primes in number
fields, while in the geometric context we will provide the desingularization of a curve.

Let $(\oo,v)$ be a valuation ring with fraction field $K$ and value group $\Gamma_v=\Z$. Let $p\in \oo$ be a uniformizing element.
Consider the following polynomials in $\oo[x]$:
\begin{equation}\label{poli}
\ars{1.2}
\begin{array}{rl}
\phi_0=&x,\\ 
\phi_1=&x^5+p^3,\\
\phi_2=&\phi_1^3+p^{10}=x^{15}+3p^3x^{10}+3p^6x^5+p^9+p^{10},   \\ 
\phi_3=&\phi_2^2+p^{11}\phi_0^4\phi_1^2, \\ 
%=&x^{30}+ 6p^3x^{25}+ 15p^6x^{20} + 2(p+10)p^{9}x^{15}+3(2p+5)p^{12}x^{10}+\\ \\
%&p^{11}x^8+2p^{14}x^6 + +6(p+1)p^{15}x^5+p^{17}x^4+(p+1)^2 p^{18}.
= & x^{30}+ 6p^3x^{25}+ 15p^6x^{20} + 2p^{9}(p+10)x^{15}+ p^{11}x^{14} +\\
&(6p+15)p^{12}x^{10}+2p^{14}x^9  +6(p+1)p^{15}x^5+p^{17}x^4+p^{18}(p+1)^2.
%=&x^{30}+ 6p^3x^{25}+ 15p^6x^{20}  + 2p^10x^{15}+ 20p^9x^{15}+ (6p+15)p^{12}x^{10}+\\\\
%&p^{11}x^8 + 2p^{14}x^6  +6(p+1)p^{15}x^5+p^{17}x^4+p^{18}(p+1)^2.
\end{array}
\end{equation}
Let us build a MLV chain of valuations on $K[x]$. We start with the  valuation:
$$
\displaystyle\mu_0\left(\sum\nolimits_i a_i x^i\right):=\min_{i} \{ v(a_i)+(3/5)i\},
$$
and consider the augmentations:
$$
\displaystyle\mu_1=[\mu_0; \phi_1,10/3], \qquad
\displaystyle\mu_2=[\mu_1; \phi_2, 301/30],  \qquad
\displaystyle\mu_3=[\mu_2; \phi_3, \infty]. 
$$

Note that  these valuations are distributed along a path of the valuative tree $\ttt(\Q)$, reaching the finite leave $\mu_3$. We get the following MLV for $\mu_3$:
\begin{equation}\label{exampleMLVchain}
%\omega_{-\infty} \stackrel{x,3/5}\lra
\mu_0 \stackrel{\phi_1,10/3}\lra \mu_1
\stackrel{\phi_2,301/30}\lra \mu_2 \stackrel{\phi_3,\infty}\lra\ \mu_3.
\end{equation}

Let us look at this MLV chain in an arithmetic context, by taking    $K=\Q$ as our base field, fixing a rational prime $p\in\Z$ and considering the $p$-adic valuation as our valuation $v$. In this setting, the MLV chain (\ref{exampleMLVchain}) encodes an important amount of information on the ring of integers $\Z_L$ of the number field $L=\Q[x]/(\phi_3)$.  For instance,  the prime ideal decomposition of $p$ in $\Z_L$ is completely described by (\ref{exampleMLVchain}). This can be checked by applying the OM-algorithm \cite{gen} to the pair $K,p$. The algorithm yields an OM-representation of $\phi_3(x)$ consisting of the unique order 3 type:
$$
\ty=(y;(x,3/5,y+1);(\phi_2,5/3,y+1); (\phi_3,1/2,y+1)),
$$
which can be seen as a computational representation of the MLV chain. It exhibits key polynomials, slopes of Newton polygons and \emph{residual polynomials}. 

Since the OM-algorithm returns a unique type, we know that 
$p\Z_L$ is divided by  a unique prime ideal $\mathfrak{P}$ whose ramification index is the product of the denominators of the slopes in $\ty$: $e(\mathfrak{P}/p)=5\cdot3\cdot2=30$. The residual degree $f(\mathfrak{P}/p)=1$ is the product of the degrees of the residual polynomials in $\ty$.

For any root $\theta\in\overline{\Q}$ of $\phi_3$, we can derive from this data the following values:
\begin{equation}\label{val-Okutsu}
\bar{v}(\theta)=3/5,\qquad \bar{v}(\phi_1(\theta))=10/3, \qquad \bar{v}(\phi_2(\theta))=301/30,
\end{equation}
where $\bar{v}$ is the unique extension of $v$ to $L$.

%PART GEOMÈTRICA...

Now, suppose that $p$ is an indeterminate, so that (\ref{poli}) can be thought as the equations of the germs at the origin $(0,0)$ of plane curves  $f(p,x)=\phi_3(x)=0$, $f_i(p,x)=\phi_i(x)=0$ with $i=1,2$,  $f_0(p,x)=p=0$, $f_{-1}(p,x)=\phi_0(x)=0$, over $\C$. Take $K=\C(p)$ as our base field, equipped with the $p$-adic valuation.
The very same OM-algorithm shows that $f$ is irreducible in $\C((p))[x]$. There is a unique finite leaf $\mu_3\in\ttt(\Q)$ with support $f\C(p)[x]$, and a MLV chain of $\mu_3$ is given in  (\ref{exampleMLVchain}). Let $\nu = 30 \mu _3$. Then, for any $\phi \in \C((p))[x]$, $\nu (\phi )$ is the intersection multiplicity between the germs of curve $f=0$ and $\phi = 0$.

The data supported by this chain contains completely analogous arithmetic information about $f$, but these data have an added geometric perspective.
Indeed, the equation (\ref{poli}) has the property that the line $p=0$ cuts the curve $f=0$ only in its %{unique}
singular point $(0,0)$, being tangent at it. In this case, the OM algorithm parallels the Newton-Puiseux algorithm for desingularization \cite[Sec. 1.2]{casas}, and the slopes and key polynomials can be reinterpreted in terms of the sequence of blow-ups involved in the desingularization process. The number of finite leaves detected by the algorithm (one in our case) is the number of points in the normalized curve lying above $(0,0)$.

Any point $P$ blown-up in the desingularization process of $f=0$ gives rise to an exceptional divisor $E_P$. For the sake of simplicity, we will use the same notation $E_P$ for any strict transform of this divisor. Any such $P$ lies either on the intersection of two exceptional divisors and $P$ is called \emph{satellite}, or it lies just on only one exceptional divisor and $P$ is called \emph{free}. 
We say that a satellite point is \emph{satellite of} the last free point preceding it.
Among all the points on $f=0$, there is a first satellite point, satellite of a free point $P_1$, which is followed by a sequence of satellite points, being $Q_1$ the last of them. Now, let  $P_2$ ($P_3$) be the second (third) free point followed by some satellite point, and let $Q_2$ ($Q_3$) be last point in the sequence of satellite points following $P_2$ ($P_3$).
The points $Q_i$ are special points in the desingularization of the curve, since they are \emph{rupture} points, that is, the exceptional divisor $E_{Q_i}$ intersects three or more other components in the pull-back of the curve.
Consider the divisorial valuation $\nu_i $ whose last centre is $Q_i$ for any $3 \geq i \geq 1$. 
It turns out that $\nu _1= 5 \mu_0$, $\nu_2 = 15 \mu _1$, $\nu_3= 30 \mu_2$ and moreover
$$
\nu (\phi_0)=18,\qquad \nu (\phi_1)=100, \qquad \nu (\phi_2)=301,
$$
which are the values appearing at (\ref{val-Okutsu}) multiplied by $\nu (f_0) = \nu (p) = 30$.
Furthermore, the germ of curve $f_{i-1}(p,x)=0$ shares with $f=0$ all its singular points and some more free simple points until $P_i$, for each $3 \geq i \geq 1$.

%En aquestes darreres frases potser caldrà retocar els subíndexos també, però no ho he fet pq no ho tenia del tot clar.
%XXXXXXXX

\subsection{Nodes of depth zero}\label{subsecDepth0}

For given $a\in K$ and $\dta\in \La\infty$, the depth-zero node $\nu=\om_{a,\dta}$ is defined as 
$$
\nu\left(\sum\nolimits_{0\le s}a_s(x-a)^s\right) = \min\left\{v(a_s)+s\dta\mid0\le s\right\}.
$$

Clearly, $\om_{a,\infty}$ is a finite leaf of $\ttt$ with support $(x-a)\kx$, while for $\dta<\infty$ the valuation $\om_{a,\dta}$ is an inner node admitting $x-a$ as a key polynomial. 

Besides these (well specified) inner nodes and finite leaves, $\ttt$ may have depth-zero infinite leaves which are the stable limit of a continuous family of augmentations of stable degree one:   
$$
\mu_0\stackrel{(\rho_i)_{i\in A}}\lra \nu,
$$
where $\mu_0$ is an inner node of depth zero. By Theorem \ref{main}, all depth-zero nodes in $\ttt$ arise from either of these two ways.
%The value group of $\om_{a,\dta}$ is: $$\g_{\om_{a,\dta}}=\begin{cases}\g,&\quad\mbox{ if }\ \dta=\infty,\\\gen{\g,\dta},&\quad\mbox{ otherwise}.\end{cases}$$

For any fixed $a \in K$, the set $\La\infty$ parametrizes a certain path in $\ttt$, containing all depth-zero nodes $\om_{a,\dta}$:

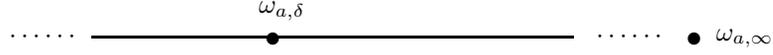
\begin{center}
\setlength{\unitlength}{4mm}
\begin{picture}(20,3)
\put(0.25,1.2){\line(1,0){16}}
\put(6,0.9){$\bullet$}\put(20,0.9){$\bullet$}
%\multiput(3,.5)(0,.25){22}{\vrule height2pt}
%\multiput(8,.9)(0,.25){9}{\vrule height2pt}
%\multiput(-.1,3)(.25,0){55}{\hbox to 2pt{\hrulefill }}
\put(-2.5,1){\begin{footnotesize}$\cdots\cdots$\end{footnotesize}}
\put(17,1){\begin{footnotesize}$\cdots\cdots$\end{footnotesize}}
\put(21,1.1){\begin{footnotesize}$\om_{a,\infty}$\end{footnotesize}}
\put(5.8,2){\begin{footnotesize}$\om_{a,\dta}$\end{footnotesize}}
\end{picture}
\end{center}

The node $\om_{a,\dta}$ is commensurable if and only if $\dta\in\gq\infty$.

The relative position of the paths corresponding to two different elements $a,b\in K$ is completely determined by the following easy observation:
\begin{equation}\label{balls}
\om_{a,\dta}\le\om_{b,\ep}\ \sii\ \min\{v(b-a),\ep\}\ge\dta.
\end{equation}

Thus,  the depth-zero paths in $\ttt$ determined by any two $a,b\in K$ coincide for all parameters $\dta\in\La$ such that 
$\dta\le v(b-a)$. 

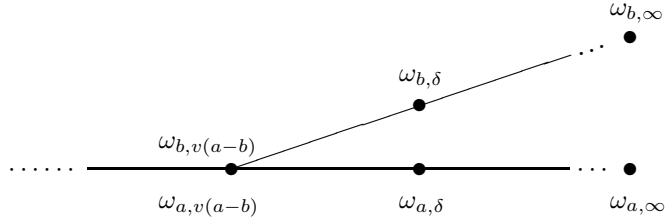
\begin{center}
\setlength{\unitlength}{4mm}
\begin{picture}(22,7.5)
\put(4.75,1){$\bullet$}\put(11,1){$\bullet$}\put(11,3.16){$\bullet$}
\put(18,1){$\bullet$}\put(18,5.4){$\bullet$}

\put(0.25,1.3){\line(1,0){16}}\put(5,1.3){\line(3,1){11.3}}
\put(16.66,4.48){$\dot{}$}\put(16.99,4.59){$\dot{}$}\put(17.32,4.7){$\dot{}$}
%\multiput(3,.5)(0,.25){22}{\vrule height2pt}
%\multiput(8,.9)(0,.25){9}{\vrule height2pt}
%\multiput(-.1,3)(.25,0){55}{\hbox to 2pt{\hrulefill }}
\put(-2.4,1.05){\begin{footnotesize}$\cdots\cdots$\end{footnotesize}}
\put(16.5,1.05){\begin{footnotesize}$\cdots$\end{footnotesize}}
\put(17.6,0){\begin{footnotesize}$\om_{a,\infty}$\end{footnotesize}}
\put(17.6,6.4){\begin{footnotesize}$\om_{b,\infty}$\end{footnotesize}}
\put(2.6,2){\begin{footnotesize}$\om_{b,v(a-b)}$\end{footnotesize}}
\put(2.6,0){\begin{footnotesize}$\om_{a,v(a-b)}$\end{footnotesize}}
\put(10.6,0){\begin{footnotesize}$\om_{a,\dta}$\end{footnotesize}}
\put(10.6,4.2){\begin{footnotesize}$\om_{b,\dta}$\end{footnotesize}}
\end{picture}
\end{center}\e

In particular, for all depth-zero nodes $\mu_0,\nu_0\in \ttt$, there is a depth-zero node $\rho_0$ such that $\rho_0<\mu_0$ and $\rho_0<\nu_0$.  By Theorem \ref{main}, for all nodes $\mu,\nu\in\ttt$ we have
\begin{equation}\label{intersect}%$$
(-\infty,\mu\,]\cap (-\infty,\nu\,]\ne\emptyset.
\end{equation}%$$

\subsection{Paths of constant depth obtained by ordinary augmentations}\label{subsecConstDepthOrd}
Let us fix an inner node $\mu\in\ttt$.
For all $\phi\in\kpm$, we define the \emph{constant-depth path} beyond $\mu$ as the set:
$$
\pmph=\left\{\mu(\phi,\ga)\mid \mu(\phi)<\ga\le\infty\right\},\qquad \mu(\phi,\ga):=[\mu;\,\phi,\ga],
$$
containing all ordinary augmentations of $\mu$ with respect to $\phi$. 
This path joins $\mu$ with the finite leaf $\mu(\phi,\infty)$. Actually, by \cite[Lem. 2.7]{MLV}, $\pmph$ coincides with the semiopen interval $(\,\mu,\mu(\phi,\infty)\,]$ in $\ttt$. 

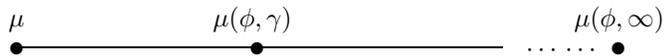
\begin{center}
\setlength{\unitlength}{4mm}
\begin{picture}(18,4)
\put(-2,1){$\bullet$}\put(-1.6,1.3){\line(1,0){16}}\put(18,1){$\bullet$}
\put(6,1){$\bullet$}
%\multiput(3,.5)(0,.25){22}{\vrule height2pt}
%\multiput(8,.9)(0,.25){9}{\vrule height2pt}
%\multiput(-.1,3)(.25,0){55}{\hbox to 2pt{\hrulefill }}
\put(-2,2){\begin{footnotesize}$\mu$\end{footnotesize}}
\put(16.5,1){\begin{footnotesize}$\cdots$\end{footnotesize}}
\put(15.2,1){\begin{footnotesize}$\cdots$\end{footnotesize}}
\put(16.8,2){\begin{footnotesize}$\mu(\phi,\infty)$\end{footnotesize}}
\put(4.8,2){\begin{footnotesize}$\mu(\phi,\ga)$\end{footnotesize}}
\end{picture}
\end{center}\e

Regardless of the commensurability or incommensurability of $\mu$, Lemma \ref{propertiesAug},(2) shows that
$\mu(\phi,\ga)$ is commensurable if and only if $\ga\in\gq\infty$.\e

\defn A key polynomial  $\phi\in\kpm$ is said to be \emph{strong} if $\deg(\phi)>\deg(\mu)$.
We say that $\pmph$ is \emph{strong} if $\phi$ is strong, and   $\pmph$  is \emph{weak} otherwise. \e

All the nodes in this path have the same degree: $\deg\left(\mu(\phi,\ga)\right)=\deg(\phi)$, and the same depth too:
$$
\dep\left(\mu(\phi,\ga)\right)=
\begin{cases}
\dep(\mu),&\mbox{ if the path is weak},\\
\dep(\mu)+1,&\mbox{ if the path is strong}.
\end{cases}
$$
Actually, for any given MLV chain of $\mu$:
$$ \mu_0\ \stackrel{\phi_1,\ga_1}\lra\  \mu_1\ \stackrel{\phi_2,\ga_2}\lra\ \cdots\ \lra\ \mu_{r-1}\ \stackrel{\phi_{r},\ga_{r}}\lra\ \mu_{r}=\mu,$$
we may obtain a MLV chain of $\mu(\phi,\ga)$ as follows. 

If the path is weak and $\mu$ has depth zero, then all $\mu(\phi,\ga)$ have depth zero too.

If  the path is weak and $\mu$ has a positive depth, we may consider
$$ \mu_0\ \stackrel{\phi_1,\ga_1}\lra\  \mu_1\ \stackrel{\phi_2,\ga_2}\lra\ \cdots\ \lra\ \mu_{r-1}\ \stackrel{\phi,\ga}\lra\ \mu(\phi,\ga),$$
regardless of the fact that $\mu_{r-1}\to\mu_r$ is an ordinary or a limit augmentation. %Therefore, the weak paths beyond $\mu$ are pieces of paths of constant depth of a lower node in $\ttt$. 

If the path is strong, we may just add one more (ordinary) augmentation: 
$$ \mu_0\ \stackrel{\phi_1,\ga_1}\lra\  \mu_1\ \stackrel{\phi_2,\ga_2}\lra\ \cdots\ \lra\ \mu_{r-1}\ \stackrel{\phi_{r},\ga_{r}}\lra\ \mu\ \stackrel{\phi,\ga}\lra\ \mu(\phi,\ga).$$

Therefore, the nodes in a strong path are ``properly" derived from $\mu$, while the nodes in weak paths are derived from lower nodes. 

%Consider the set of all nodes in {\bf strong} paths of constant depth  beyond $\mu$:$$\pmu=\bigcup_{\phi\in\kpmz}\pmph\ \subset\  \ttt.$$Note that $\mu\not\in\pmu$. 

Let us analyze the intersection of two paths of constant depth beyond the same node $\mu$, determined by different key polynomials $\phi,\phi_*\in\kpm$.
Obviously, 
$$\ty(\mu,\rho)=[\phi]_\mu\quad \mbox{for all}\quad \rho\in \pmph.
$$
Therefore, if $\phi\not\sim_\mu\phi_*$, Proposition \ref{td=td} shows that $\pmph\cap\pset_\mu(\phi_*)=\emptyset$.

If $\phi\smu\phi_*$, then $\mu(\phi-\phi_*)>\mu(\phi)=\mu(\phi_*)$. Thus, (\ref{eqAug}) shows that 
$$\pmph\cap\pset_\mu(\phi_*)=\left(\mu,\mu(\phi,\ga_0)\right],\qquad \ga_0=\mu(\phi-\phi_*).
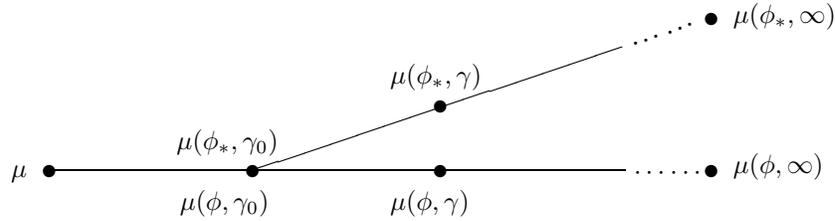
$$

\begin{center}
\setlength{\unitlength}{4mm}
\begin{picture}(22,8)
\put(-2,1){$\bullet$}\put(4.75,1){$\bullet$}
\put(11,1){$\bullet$}\put(11,3.16){$\bullet$}
\put(20,1){$\bullet$}\put(20,6.06){$\bullet$}

\put(-1.6,1.3){\line(1,0){19}}\put(5,1.3){\line(3,1){12.3}}
\put(17.76,4.88){$\dot{}$}\put(18.09,4.99){$\dot{}$}\put(18.42,5.11){$\dot{}$}
\put(18.9,5.26){$\dot{}$}\put(19.23,5.37){$\dot{}$}\put(19.56,5.5){$\dot{}$}
%\multiput(3,.5)(0,.25){22}{\vrule height2pt}
%\multiput(8,.9)(0,.25){9}{\vrule height2pt}
%\multiput(-.1,3)(.25,0){55}{\hbox to 2pt{\hrulefill }}
\put(-3,1.1){\begin{footnotesize}$\mu$\end{footnotesize}}
\put(17.7,1){\begin{footnotesize}$\cdots$\end{footnotesize}}
\put(18.9,1){\begin{footnotesize}$\cdots$\end{footnotesize}}
\put(21,1.2){\begin{footnotesize}$\mu(\phi,\infty)$\end{footnotesize}}
\put(21,6.2){\begin{footnotesize}$\mu(\phi_*,\infty)$\end{footnotesize}}
\put(2.5,2){\begin{footnotesize}$\mu(\phi_*,\ga_0)$\end{footnotesize}}
\put(2.6,0){\begin{footnotesize}$\mu(\phi,\ga_0)$\end{footnotesize}}
\put(9.6,0){\begin{footnotesize}$\mu(\phi,\ga)$\end{footnotesize}}
\put(9.6,4.2){\begin{footnotesize}$\mu(\phi_*,\ga)$\end{footnotesize}}
\end{picture}
\end{center}\bs

\subsection{Paths of constant depth obtained by limit augmentations}\label{subsecConstDepthLim}
Let us fix an essential continuous family $\aa=\left(\rho_i\right)_{i\in A}$. By taking a cofinal family, if necessary, we may assume that $\aa$ contains a minimal valuation $\mu$.

For all limit key polynomials $\phi\in\kpi(\aa)$, we may consider the \emph{constant depth path} beyond $\aa$:
$$
\paph=\left\{\aa(\phi,\ga)\mid \rho_i(\phi)<\ga\le\infty \ \mbox{ for all }i\in A\right\},\qquad \aa(\phi,\ga):=[\aa;\,\phi,\ga],
$$
containing all possible limit augmentations determined by $\phi$.

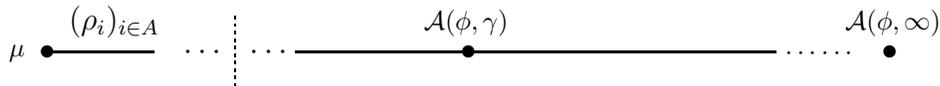
\begin{center}
\setlength{\unitlength}{4mm}
\begin{picture}(26,4)
\put(-2,1){$\bullet$}\put(26,1){$\bullet$}\put(12,1){$\bullet$}
\put(-1.6,1.25){\line(1,0){3.4}}\put(2.8,1){$\cdots$}\put(5,1){$\cdots$}\put(6.5,1.25){\line(1,0){16}}
\multiput(4.5,0.1)(0,.25){10}{\vrule height1pt}
\put(-1,2){$(\rho_i)_{i\in A}$}
%\multiput(8,.9)(0,.25){9}{\vrule height2pt}
%\multiput(-.1,3)(.25,0){55}{\hbox to 2pt{\hrulefill }}
\put(-3,1.1){\begin{footnotesize}$\mu$\end{footnotesize}}
\put(22.8,1){\begin{footnotesize}$\cdots$\end{footnotesize}}
\put(24,1){\begin{footnotesize}$\cdots$\end{footnotesize}}
\put(24.8,2){\begin{footnotesize}$\aa(\phi,\infty)$\end{footnotesize}}
\put(10.8,2){\begin{footnotesize}$\aa(\phi,\ga)$\end{footnotesize}}
\end{picture}
\end{center}\e

As in the previous cases, the last node of the path, $\aa(\phi,\infty)$, is a finite leaf. 
By Proposition \ref{extensionlim}, $\aa(\phi,\ga)$ is commensurable if and only if $\ga\in\gq\infty$.
By \cite[Lem. 3.8]{MLV}, we have $\paph=\bigcap\nolimits_{i\in A}(\,\rho_i,\aa(\phi,\infty)\,]$.

All the nodes in this path have the same degree and the same depth: 
$$\deg\left(\aa(\phi,\ga)\right)=\deg(\phi),\qquad\dep\left(\aa(\phi,\ga)\right)=\dep(\mu)+1.
$$

%Consider the set of all nodes in all paths of constant depth  beyond $\aa$:$$\paa=\bigcup_{\phi\in\kpi(\aa)}\paph \subset\ \ttt,$$

For any $\phi,\phi_*\in\kpi(\aa)$, let $\dta=\rho_\aa(\phi-\phi_*)$. By \cite[Lem. 3.7]{MLV}, we have
$$
\aa(\phi,\ga)=\aa(\phi_*,\ga_*)\sii \ga=\ga_*\ge\dta.
$$
By Lemma \ref{allLKP}, $\aa(\phi,\dta)=\aa(\phi_*,\dta)$ belongs to $\paph\cap\pset_\aa(\phi_*)$.

%The following result, which says  in some sense that $\aa$ has a unique tangent direction, is an immediate consequence of Lemma \ref{allLKP}.

%\begin{lemma}\label{OneTd}For all $\phi,\phi_*\in\kpi(\aa)$, we have $$\rho_\aa(\phi-\phi_*)>\rho_i(\phi)=\rho_i(\phi_*)\quad \mbox{ for all }\ i\in A.$$\end{lemma}

Therefore, the intersection of the paths determined by $\phi$ and $\phi_*$ is completely analogous to the case of depth-zero valuations.

\subsection{Greatest common lower node}\label{subsecGCN}

Given  $\mu,\nu\in\ttt$, their \emph{greatest common lower node} is defined as
$$
\mu\wedge\nu=\max\left( (-\infty,\mu\,]\cap (-\infty,\nu\,]\right)\in \ttt,
$$
provided that this maximal element exists.

\begin{proposition}\label{GCN}
For all $\mu,\nu\in\ttt$, their greatest common lower node $\mu\wedge\nu$ exists
\end{proposition}

\begin{proof}
If $\mu\le\nu$, then obviously $\mu\wedge\nu=\mu$. Suppose that neither $\mu\le\nu$ nor $\mu\ge\nu$. 
As we saw in (\ref{intersect}), 
$
(-\infty,\mu\,]\cap (-\infty,\nu\,]\ne\emptyset$.
Let us prove that this totally ordered set always contains a maximal element. 

Consider a MLV chain of $\mu$ %and $\nu$. 
$$
\mu_0\ \stackrel{\phi_{1},\ga_{1}}\lra\  \mu_{1}\ \lra\ \cdots
\ \lra\ \mu_{r-1} 
\ \stackrel{\phi_{r},\ga_{r}}\lra\ \mu_{r}=\mu, 
$$
%$$\mu_0^*\ \stackrel{\phi^*_{1},\ga^*_{1}}\lra\  \nu_{1}\ \lra\ \cdots\ \lra\ \nu_{s-1} \ \stackrel{\phi^*_{s},\ga^*_{s}}\lra\ \nu_{s}=\nu.$$
Since $\mu=\mu_r\not\in (-\infty,\nu\,]$, there exists a minimal index $i$ such that $\mu_i\not\in (-\infty,\nu\,]$.
We need only to show that $\mu_i\wedge\nu$ exists, because this clearly implies $\mu\wedge\nu=\mu_i\wedge\nu$.

Suppose that $i=0$. If $\mu_0=\om_{a,\ga}$, we have \cite[Sec. 2.2]{MLV}
$$
\left(-\infty,\mu_0\,\right]=\left\{\om_{a,\dta}\mid\dta\in\La,\ \dta\le\ga\right\}.
$$
On the other hand, by comparing their action of $(x-a)$-expansions, we see that $\om_{a,\dta}\le\nu$ if and only if $\dta\le\nu(x-a)$. Since $\mu_0\not\le\nu$, necessarily $\nu(x-a)<\ga$, Thus, there is a maximal element in $\left(-\infty,\mu_0\,\right]\cap\left(-\infty,\nu\,\right]$, namely
$$\mu_0\wedge\nu=\om_{a,\nu(x-a)}.$$

Suppose that $i>0$, so that $\mu_{i-1}<\nu$, $\mu_{i}\not\le\nu$. 

If $\ty(\mu_{i-1},\mu_i)\ne\ty(\mu_{i-1},\nu)$, then Proposition \ref{td=td} shows that $(\mu_{i-1},\mu_i\,]\cap (\mu_{i-1},\nu]=\emptyset$. Hence, $\mu_i\wedge\nu=\mu_{i-1}$. Suppose that $\ty(\mu_{i-1},\mu_i)=\ty(\mu_{i-1},\nu)$.

If $\mu_{i-1}\to\mu_i$ is an ordinary augmentation, then \cite[Lem. 2.7]{MLV} shows that
$$
\left(\mu_{i-1},\mu_i\,\right]=\left\{[\mu_{i-1};\,\phi_i,\dta]\,\mid\,\mu_{i-1}(\phi_i)<\dta\le\ga_i\right\}.
$$
Since $\ty(\mu_{i-1},\nu)=\ty(\mu_{i-1},\mu_i)=[\phi_i]_{\mu_{i-1}}$, we have $\mu_{i-1}(\phi_i)<\nu(\phi_i)$ and
$$
\mu_{i-1}(f)=\nu(f) \sii\phi_i\nmid_{\mu_{i-1}}f.
$$
In particular, $\mu_{i-1}(a)=\nu(a)$ for all $a\in\kx$ with $\deg(a)<\deg(\phi_i)$. Hence, by comparing their action on $\phi_i$-expansions, we have
$$
[\mu_{i-1};\,\phi_i,\dta]\le\nu\sii \dta\le \nu(\phi_i).
$$
Since $\mu_i\not\le\nu$, we have $\nu(\phi_i)<\ga_i$. Thus,  there is a maximal element in $\left(\mu_{i-1},\mu_i\,\right]\cap\left(\mu_{i-1},\nu\,\right]$, namely
$$\mu_i\wedge\nu=[\mu_{i-1};\,\phi_i,\nu(\phi_i)].$$

Suppose that $\mu_{i-1}\to\mu_i$ is a limit augmentation with respect to  an  essential continuous family $\aa=(\rho_j)_{j\in A}$ admitting $\mu_{i-1}$ as its first element. Then, $\phi_i\in\kpi(\aa)$. By  Lemma \ref{specialCont} we may assume that, for all $j\in A$, 
$$
\rho_j=[\mu_{i-1};\,\chi_j,\be_j],\qquad \chi_j\in\kp(\mu_{i-1}),\ \be_j=\sval(\rho_j).
$$
%We distinguish two cases, according to the existence, or not,  of some $j\in A$ such that $\rho_j\not\le\nu$.

If $\rho_j\not\le\nu$ for some $j\in A$, then we can mimic the arguments of the ordinary-augmentation case to conclude that 
$$\mu\wedge\nu=\rho_j\wedge\nu=[\mu_{i-1};\,\chi_j,\nu(\chi_j)].$$ 

Suppose that $\rho_j<\nu$ for all $j\in A$. By Lemma \ref{propertiesTMN}(3), we see that  
 $\nu$ coincides with $\rha$ on all $\aa$-stable polynomials.
Let $V=\left\{\rho_j(\phi_i)\mid j \in A\right\}$. By \cite[Lem. 3.8]{MLV},
$$
\left\{\rho\in(\mu_{i-1},\mu_i\,]\,\mid\,\rho_j<\rho\mbox{ for all }j\in A\right\}=
\left\{[\aa;\,\phi_i,\dta]\mid V<\dta\le\ga_i\right\}.
$$
By comparing their action on $\phi_i$-expansions, we have
$$
[\aa;\,\phi_i,\dta]\le\nu\sii \dta\le\nu(\phi_i).
$$
On the other hand, for all $\aa$-unstable polynomials, we have  $\rho_j(f)<\rho_\ell(f)$ for all $j<\ell$ in $A$. Thus, 
$$\rho_j(f)<\rho_\ell(f)\le\nu(f),\quad\mbox{for all }j\in A. $$
In particular, $\nu(\phi_i)>V$.
As a consequence,  there is a maximal valuation in $(-\infty,\mu_i\,]$ which is less than $\nu$, namely
$$
\mu_i\wedge\nu=[\aa;\,\phi_i,\nu(\phi_i)].
$$
This ends the proof of the proposition.
\end{proof}\e

Suppose that $\mu,\nu\in\ttt$ are incomparable; that is, $\mu\not\le\nu$ and $\nu\not\le\mu$. Then,
 their greatest common lower node, $\rho=\mu\wedge\nu$, has at least two different tangent directions: $\ty(\rho,\mu)\ne \ty(\rho,\nu)$. By Theorems \ref{Delta} and \ref{DeltaKP}, $\rho$ is an inner commensurable node; in other words, $\rho$ is a residually transcendental valuation. 

\subsubsection{$\ttt$ as a $\Lambda$-tree}

Given an ordered abelian group $\Lambda$, a $\Lambda$-tree is defined \cite{lambdatrees} as a geodesic $\Lambda$-metric space $T$ such that 
\begin{enumerate}
    \item If two geodesics of $T$ intersect in a single point, which is an endpoint of both, then their union is a geodesic;
    \item The intersection of two geodesics with a commond endpoint is also a geodesic.
\end{enumerate}
The existence of a greatest common lower node can be used to define a $\Lambda$-metric on the subtree $\mathcal{H}$ of $\ttt$ consisting of all inner nodes.
Namely, we set
\[
d(\mu,\nu)=
\sval(\mu)+\sval(\nu)-2\,\sval(\mu\wedge\nu). 
\]
Note that $d(\mu,\nu)=\left|\sval(\mu) - \sval(\nu)\right|$ if $\mu$ and $\nu$ are comparable.
It is easy to see that with this definition, $\mathcal{H}$ is a geodesic $\Lambda$-metric space, and the unique geodesic with endpoints $\mu$, $\nu$ is the union of the segments $[\mu\wedge\nu,\mu]$ and $[\mu\wedge\nu,\nu]$; the two properties above follow.

We are not going to use any metric properties of the tree $\mathcal{H}$, noting only that this is a hyperbolic space. This fact, along with a plethora of additional information, can be found in the monograph \cite{lambdatrees}.

\section{Equivalence classes of valuations and small extensions of groups}\label{secSME}

For our given valued field $(K,v)$, consider a valuation $\mu\colon \kx\to\La\infty$,
whose restriction to $K$ is equivalent to $v$. That is, there exists an order-preserving embedding $j\colon \g\hk\La$, fitting into a commutative diagram
$$
\ars{1.3}
\begin{array}{ccc}
\kx&\stackrel{\mu}\lra&\La\infty\\
\uparrow&&\ \uparrow\mbox{\tiny$j$}\\
K&\stackrel{v}\lra&\g\infty
\end{array}
$$

The induced embedding  $j\colon \g\hk\gm$ is necessarily a \emph{small extension} of ordered abelian groups. That is, if $\g'\subset\gm$ is the relative divisible closure of $\g$ in $\gm$, then $\gm/\g'$ is a cyclic group \cite[Thm. 1.5]{Kuhl}.

Not all small extensions of $\g$ arise from valuations on a polynomial ring. In \cite{Kuhl} it is shown that the divisible closure of $\g$ in $\gm$ must be countably generated over $\g$, and it must be finitely generated over $\g$, if $\gm/\g$ is not a torsion group.  

%This theorem was proved in  \cite[Thm. 1.5]{Kuhl} for valuations with  trivial support; that is, valuations that may be extended to the rational field $K(x)$. For  valuations with non-trivial support the extension $\g\subset\gm$ is commensurable, because $\gm$ is the value group of an extension of $\mu_{\mid K}$ to a finite extension of $K$.

Our aim is to describe the tree $\ttt_v$ whose nodes are all equivalence classes of va\-lua\-tions on $\kx$ whose restriction to $K$ is equivalent to $v$. The first natural step is to build up some universal ordered group $\La$ containing all small extensions of $\g$ up to order-preserving $\g$-isomorphism.

%\subsubsection*{Principal convex subgroups of $\g$ and Hahn's embedding theorem}\mbox{\null}\e

\subsection{Maximal rank-preserving extension of $\g$}\label{subsecRlex}

From now on, an \emph{embedding} of totally ordered sets is a mapping which strictly preserves the ordering. Also, an \emph{embedding} $\La\hk \La'$ of totally ordered abelian groups is a group homomorphism which is an embedding as totally ordered sets.

A subgroup $H\subset \g$ is \emph{convex} if for all positive $h\in H$, it holds $[-h,h]\subset H$. For all $a\in \g$, the intersection of all convex subgroups of $\g$ containing $a$, is a \emph{principal} convex subgroup of $\g$.

Let $\pcv(\g)$ be the totally ordered set of {\bf nonzero} principal convex subgroups of $\g$, ordered by {\bf decreasing} inclusion.

Any embedding $j\colon \g\hk\La$ induces an embedding of ordered sets$$\pcv(\g)\hk\pcv(\La),$$
which maps the principal convex subgroup generated by $a\in\g$ to the principal convex subgroup generated by $j(a)$ in $\La$.\e

\defn
We say that $j$ \emph{preserves the rank} if this mapping is bijective. \e

For instance, the canonical embeddding $\g\hk\gq$ preserves the rank. From now on, we shall consider the bijection between $\pcv(\g)$ and $\pcv(\gq)$ as an identity:
$$
I:=\pcv(\g)=\pcv(\gq).
$$

We may identify $I\infty$ with a set of indices parameterizing all principal convex subgroups of $\gq$. For all $i\in I$ we denote by $H_i$ the corresponding  principal convex subgroup. We agree that $H_\infty=\{0\}$.
Then, according to our convention, for any pair of indices $i,j\in I\infty$, we have
$$
i<j \sii H_i\supsetneq H_j.
$$

%The order-type of $I$ is called the \emph{principal rank} of $\g$, and is denoted $\pr(\g)$.\e

Let $\left\{I,(Q_i)_{i\in I}\right\}$ be the skeleton of $\gq$. That is, $Q_i=H_i/H_i^*$ for all $i\in I$, where $H_i^*\subset H_i$ is the maximal proper convex subgroup of $H_i$. That is, if $a\in H_i$ generates $H_i$ as a convex subgroup, then $H_i^*$ is the union of all convex subgroups of $\gq$ not containing $a$.
The convex subgroup $H_i^*$ is not necessarily principal.

Consider the respective Hahn's products:
$$
\hq\subset \prod\nolimits_{i \in I}Q_i,\qquad \rlx\subset \R^I,
$$
equipped with the lexicographical order. That is, $\hq,\,\rlx$ are the subsets of $\prod_{i \in I}Q_i,\,\R^I$, respectively, formed by all elements $x=(x_i)_{i\in I}$ whose support
$$
\supp(x)=\{i\in I\mid x_i\ne0\}
$$is a well-ordered subset of $I$, with respect to the ordering induced by that of $I$.

By Hahn's embedding theorem \cite[Sec. A]{Rib}, there exists a (non-canonical) $\Q$-linear embedding $$\gq\hra\hq$$ which induces an isomorphism between the respective skeletons.

On the other hand, the ordered groups $Q_i$ have rank one; that is, they have only two convex subgroups: $\{0\}$ and $Q_i$.
Hence, the choice of positive elements $1_i\in Q_i$ determines $\Q$-linear embeddings for all $i\in I$:
$$%\begin{equation}\label{QiR}
Q_i\hra \R,\qquad 1_i\longmapsto 1,
$$%\end{equation}

Therefore, we have a natural embedding $\hq\hk \rlx$.
All in all, we obtain a rank-preserving extension
$$
\tau\colon \g\hra\gq\hra\hq\hra\rlx,
$$
which is maximal among all rank-preserving extensions of $\g$ \cite[Sec. A]{Rib}.

\begin{theorem}\label{MaxEqRk}
For any rank-preserving extension $\g\hk \La$, there exists an embedding $\La\hk\rlx$ fitting into a commutative diagram
$$
\ars{1.2}
\begin{array}{ccc}
\La&&\\
\uparrow&\raise.5ex\hbox{$\searrow$}&\\
\g&\stackrel{\tau}\lra&\rlx
\end{array}
$$
\end{theorem}

The embedding $\La\hk\rlx$ is not unique. Thus, every rank-preserving extension of $\g$ is $\g$-equivalent to some subgroup of $\rlx$, but not to a unique one.

The nonzero principal convex subgroups of $\rlx$ are parametrized by $I$ via:
$$
H_i=\{(x_j)_{j\in I}\mid x_j=0\mbox{ for all }j<i\}\quad\mbox{for all }i\in I.
$$

The convex subgroups of $\rlx$ are parametrized by the set $\inii$ of initial segments of $I$, as follows:
\begin{equation}\label{HSHahn}
S\in \inii\quad \rightsquigarrow\quad H_S=\{(x_j)_{j\in I}\mid x_j=0\mbox{ for all }j\in S\}.
\end{equation}

\subsection{A universal group for small extensions of $\g$}\label{subsecRii}

For all $S\in\inii$, let $i_S$ be a formal symbol and consider the ordered set
$$
I_S=S+\{i_S\}+(I\setminus S),
$$
where $+$ is the usual addition of totally ordered sets.

We define the \emph{one-added-element hull} of $I$ as the set
$$
\I:=I\cup\left\{i_S\mid S\in\inii\right\},
$$
equipped with the total ordering determined by
\begin{enumerate}
\item[(i)] For all $S\in\inii$, the inclusion $I_S\hk \I$ preserves the order.
\item[(ii)] For all $S,T\in\inii$, we have $i_S<i_T$ if and only if $S\subsetneq T$.
\end{enumerate}

Consider the Hahn product $\rii\subset \R^{\I}$, defined as above, just by replacing $I$ with $\I$.

The inclusions $I\subset  I_S\subset \I$ determine canonical embeddings
$$
\rlx\hra \rlxs\hra \rii,\quad\mbox{for all }S\in\inii.
$$
Altogether, we obtain an embedding
$$
\ell\colon \g\hra\gq\hra\hq\hra\rlx\hra\rii.
$$

As shown in \cite[Prop. 5.1]{csme}, $\rii$ is the universal ordered abelian group we are looking for.

\begin{proposition}\label{riiUniverse}
Let $\mu$ be a valuation on $\kx$ whose restriction to $K$ is equivalent to $v$. Then, there exists an embedding $j\colon\gm\hk \rii$ satisfying the following properties:
\begin{enumerate}
\item[(i)] The following diagram commutes:
$$
\ars{1.3}
\begin{array}{ccc}
\kx&\stackrel{j\circ\mu}\lra&\rii\infty\\
\uparrow&&\ \uparrow\mbox{\tiny$\ell$}\\
K&\stackrel{v}\lra&\g\infty
\end{array}
$$
\item[(ii)] There exists $S\in\inii$ such that $j(\gm)\subset \rlxs$.
\end{enumerate}

Moreover, $\gm/\g$ is commensurable if and only if $j(\gm)\subset\ell(\g)_\Q$. Also, $\gm/\g$ preserves the rank if and only if $j(\gm)\subset\rlx$.
\end{proposition}

Clearly, $\mu$ is equivalent to the valuation $j\circ \mu$ on $\kx$, and $v$ is equivalent to the valuation $\ell\circ v$ on $K$. Also, the valuation $j\circ \mu$ restricted to $K$ is equal to $\ell\circ v$.

As a consequence, in order to describe all equivalence classes of valuations $\mu$ on $\kx$ whose restriction to $K$ is equivalent to $v$, we may assume that $v$ and $\mu$ satisfy the following conditions:\e

(V1) \ The valuation $v$ takes values in $\rlx$. That is, $\g=v(K^*)\subset \rlx$.\e%\footnote{In other words, the inclusions $\g\subset \gq\subset \hq\subset \rlx$ concern adequate subgroups of $\rlx$.}\e

(V2) \ The valuation $\mu$ satisfies $\mu_{\mid K}=v$ and takes values in $\rlxs$ for some $S\in\inii$.%\e

%To this end will be devoted Section \ref{secConstruct}.

\subsection{Small-extensions equivalence on a subset of the universal value group}\label{subsecSME}
From now on, we assume that our valuation $v$ on $K$ satisfies  (V1), so that the embedding $\ell$ of the last section is the canonical inclusion.
Consider the subset
$$\rll=\rll(I):=\bigcup_{S\in\inii}\rlxs\subset \rii.$$

For all $\be\in\rll$, we denote by $\ggb\subset \rii$ the subgroup generated by $\g$ and $\be$. %Clearly, if $\be\in \rlxs$, then $\ggb\subset \rlxs\subset\rll$.

Let $\La=\rii$ and $\ttt=\tla$. Consider the subtree
$$\tz=\left\{\rho\in\ttt\mid \grh\subset\rlxs\ \mbox{ for some }\ S\in\inii \right\}\subset \ttt.$$
Note that all valuations in $\tz$ satisfy the condition (V2).

On the set $\rll$ we define the following equivalence relation.\e

\noindent{\bf Definition.} We say that $\be,\ga\in \rll$ are $\,\op{sme}$-equivalent if there exists an isomorphism of ordered groups
$$
\ggb\iso\gga,
$$
which acts as the identity on $\g$ and sends $\be$ to $\ga$.

In this case, we write $\be\sim_\sme\ga$. %if the base group $\g$ is clear from the context.
We denote by $[\be]_\sme\subset \rll$ the class of $\be$.\e

The motivation for this definition lies in the following result.

\begin{proposition}\label{motivation}
Let $\mu,\nu\in\tz$ be two inner nodes. Then, $\mu\sim\nu$ if and only if the following three conditions hold:

(a) \ The valuations $\mu$, $\nu$ admit a common key polynomial of minimal degree.

(b) \ For all \,$a\in\kx\,$ such that \,$\deg(a)<\deg(\mu)$, we have $\,\mu(a)=\nu(a)$.

(c) \ $\sval(\mu)\sim_\sme \sval(\nu)$.

In this case, we have $\kpm=\kpn$.
\end{proposition}

\begin{proof}
Suppose that $\mu\sim\nu$. Then, there exists an isomorphism of ordered groups $\iota\colon \gm\ism\gn$ such that $\nu=\iota\circ\mu$.
The isomorphism $\iota$ induces an isomorphim between the graded algebras:
$$
\ggm\iso\ggn,\qquad f+\pset^+_\al(\mu)\longmapsto f+\pset^+_{\iota(\al)}(\nu),
$$
for all $f\in\pset_\al(\mu)$, $\al\in\gm$. Since key polynomials are characterized by algebraic properties of their initial terms in the graded algebra, this implies that both valuations have the same key polynomials: \,$\kpm=\kpn$.

Let $\phi$ be a common key polynomial of minimal degree.
Since $\mu_{\mid K}=\nu_{\mid K}$, the isomorphism $\iota$ restricted to $\g$ is the identity and $\iota(\mu(\phi))=\nu(\phi)$. Hence,
$$\sval(\mu)=\mu(\phi)\sim_\sme\nu(\phi)=\sval(\nu).
$$

Finally, since  $\iota$ restricted to $\g$ is the identity, then $\iota$ acts as the identity on any torsion element in $\gm/\g$. Now, for all $a\in \kx$ of degree less than $\deg(\phi)$, the values $\mu(a)$, $\nu(a)$ belong to $\gq$ \cite[Lem. 1.3]{MLV}. Thus,
$$
\nu(a)=\iota(\mu(a))=\mu(a).
$$

Conversely, suppose that $\mu$ and $\nu$ satisfy conditions (a), (b) and (c). Take $\phi\in\kpm\cap\kpn$ of minimal degree in both sets. Let us denote $$\be=\sval(\mu)=\mu(\phi), \qquad \ga=\sval(\nu)=\nu(\phi).$$ By condition (c), there is an order-preserving $\g$-isomor\-phism $\iota\colon\gen{\g,\be}\ism\gen{\g,\ga}$, mapping $\be$ to $\ga$.
As we saw in Section \ref{subsubsecG0}, the subgroup
$$
H:=\gm^0=\left\{\mu(a)\mid 0\le \deg(a)<\deg(\mu)\right\}\subset\gm
$$
is commensurable over $\g$ and satisfies $\gm=\gen{H,\mu(\phi)}=\gen{H,\be}$. By condition (b), $H=\gn^0$ and $\gn=\gen{H,\nu(\phi)}=\gen{H,\ga}$ too. Since $H/\g$ is a torsion abelian group, the $\g$-isomorphism $\iota$ induces an order-preserving isomorphism
$$\iota\colon\gm=\gen{H,\be}\iso\gn=\gen{H,\ga},$$ which acts as the identity on $H$ and maps $\be$ to $\ga$. Let us check that $\nu=\iota\circ\mu$.

For $f\in\kx$, consider its $\phi$-expansion $f=\sum_{0\le s}a_s\phi^s$, where $\deg(a_s)<\deg(\phi)=\deg(\mu)$. Since $\phi$ is $\mu$-minimal and $\nu$-minimal, Lemma \ref{minimal0} shows that
$$
\mu(f)=\min\left\{\mu(a_s\phi^s)\mid 0\le s\right\},\qquad \nu(f)=\min\left\{\nu(a_s\phi^s)\mid 0\le s\right\}.
$$
Let us denote $\dta_s=\mu(a_s)=\nu(a_s)\in H$, for all $s\ge0$. Clearly,
$$
\iota\left(\mu(a_s\phi^s)\right)=\iota(\dta_s+s\be)=\dta_s+s\ga=\nu(a_s\phi^s).
$$
Since $\iota$ preserves the ordering, for arbitrary indices $s,t$ we have
$$
\mu(a_s\phi^s)\le\mu(a_t\phi^t)\imp \nu(a_s\phi^s)\le\nu(a_t\phi^t).
$$
Thus, there is a common index $s$ for which $\mu(f)=\mu(a_{s}\phi^{s})$ and $\nu(f)=\nu(a_{s}\phi^{s})$, simultaneously. Therefore,
$
\iota\left(\mu(f)\right)=\iota\left(\mu(a_s\phi^s)\right)=\nu(a_s\phi^s)=\nu(f)
$.
\end{proof}

%The following result shows concrete instances where the conditions of Proposition \ref{motivation} hold.

\begin{corollary}\label{motivation2}
Take $\be,\,\ga\in \rll$.
\begin{enumerate}
\item[(i)] For all $a\in K$ we have  $ \om_{a,\be}\sim\om_{a,\ga}$ if and only if $\be\sim_\sme\ga$.
%where $\om_{a,\be}$ is the depth-zero valuation introduced in Section \ref{subsecDepth0}.
\item[(ii)] Let $\mu\in\tz$ be an inner node and let $\phi\in\kpm$. If $\be,\ga>\mu(\phi)$, then
$$
[\mu;\,\phi,\be]\sim [\mu;\,\phi,\ga]\sii \be\sim_\sme\ga.
$$
\item[(iii)] Let $\aa=\left(\rho_i\right)_{i\in A}$ be an essential continuous family in $\tz$ and let $\phi\in\kpi(\aa)$. If $\be,\ga>\rho_i(\phi)$ for all $i\in A$, then
$$
[\aa;\,\phi,\be]\sim [\aa;\,\phi,\ga]\sii \be\sim_\sme\ga.
$$
\end{enumerate}
\end{corollary}

\begin{proof}
All items follow immediately from Proposition \ref{motivation}, once we see that for the two involved valuations, conditions (a) and (b) hold in each case.

In case (i), the common key polynomial of minimal degree is $\phi=x-a$.

In cases (ii) and (iii), $\phi$ is a common key polynomial of minimal degree for both augmentations  by Lemma \ref{propertiesAug} and Proposition \ref{extensionlim}, respectively. 
\end{proof}%\e

\subsection{Quasi-cuts in $\gq$ and small-extensions closure of $\g$}\label{subsecQcuts}

Consider any subset $\gsme\subset \rll$ which is a set of representatives of the quotient set $\rll/\!\sim\sme$.

The only $\g$-automorphism of $\gq$ as an ordered group is the identity. Thus, for all $\be\in\gq\subset\rll$, we have $\left[\be\right]_\sme=\{\be\}$. Therefore, we have necessarily
$$\gq\subset\gsme\subset \rll.$$

Any such ``small-extensions closure" $\gsme$ contains generators of all small extensions of $\g$, up to the relative divisible closure of $\g$. Let us give a precise explanation of this statement, which follows easily from the definition of $\sim_\sme$.

\begin{proposition}\label{smallSme}Let $\g\hk\Omega$ be a small extension and let $\ga\in\Omega$ such that $\Omega=\gen{\d,\ga}$, where $\d$ is the relative divisible closure of $\g$ in $\Omega$. Let $\d\!\ism\d_0\subset\gq$ be the canonical embedding of $\d$ into $\gq$. Then, for a  unique $\be\in\gsme$ there exists an isomorphism  of ordered groups $$\Omega\iso\gen{\d_0,\be},\qquad \ga\longmapsto \be,$$ which maps $\ga$ to $\be$, and whose restriction to $\d$ is the canonical isomorphism $\d\!\ism \d_0$.
\end{proposition}

In \cite{csme}, a real model for the set of quasi-cuts of $\gq$ is constructed, which serves as a canonical choice for $\gsme$. Let us recall this construction.

A \emph{quasi-cut} in $\gq$ is a pair $D=\left(D^L,D^R\right)$ of subsets such that $D^L\le D^R$ and $D^L\cup D^R=\gq$.
Then, $D^L$ is an initial segment of $\gq$, $D^R$ is a final segment of $\gq$ and $D^L\cap D^R$ consists of at most one element.

If $D^L\cap D^R=\{a\}$, we say that $D$ is the \emph{principal quasi-cut} determined by $a\in \gq$. If  $D^L\cap D^R=\emptyset$, we say that $D$ is a \emph{cut} in $\gq$.

The set $\op{Qcuts}(\gq)$ of all quasi-cuts in $\gq$ admits a total ordering:
$$
D=\left(D^L,D^R\right)\le E=\left(E^L,E^R\right) \ \sii\ D^L\subset E^L\quad\mbox{and}\quad D^R\supset E^R.
$$

For all $x\in\rll$ we consider the folllowing quasi-cut $D_x$ in $\gq$:
$$
D_x^L=\{a\in\gq\mid a\le x\},\qquad D_x^R=\{a\in\gq\mid a\ge x\}.
$$
We say that $x$ \emph{realizes} the cut $D_x$. The set $\rll$ contains realizations of all quasi-cuts in $\gq$ \cite[Sec. 4]{csme}.
Moreover, these quasi-cuts provide the following  reinterpretation of the equivalence relation $\sim_\sme$ \cite[Lem. 5.4]{csme}.

\begin{lemma}
For all $x,y\in\rll$, we have $x\sim_\sme y$ if and only if  $D_x=D_y$.
\end{lemma}

As a consequence, if we consider on $\gsme$ the total ordering induced by that of $\rll$,
we derive a natural isomorphism of ordered sets:
$$\gsme\lra \op{Qcuts}(\gq),\qquad x\longmapsto D_x.$$

\begin{corollary}\label{complete}
	Equipped with the order topology,  $\gsme$ is complete and contains $\gq$ as a dense subset.
\end{corollary}

Indeed, it is well known that the ordered set $\op{Qcuts}(\gq)$ has these properties. We recall that being complete with respect to the order topology means  that every non-empty subset of $\gsme$ has a supremum and an infimum.

In \cite[Sec. 4]{csme}, a canonical choice for $\gsme$ is described as
$$
\gsme=\gq\sqcup \gnbc\sqcup \gbc,
$$
for certain subsets $\gnbc\subset \rlx\setminus\gq$ and $\gbc\subset \rll\setminus\rlx$.

The elements $x\in \gq$ parametrize the principal quasi-cuts. The elements $x\in\gnbc$, $x\in\gbc$ correspond to $D_x$ being a \emph{non-ball cut}, or a \emph{ball cut}. Equivalently, the small extension $\g\hk\gen{\gq,x}$ preserves, or increases the rank, respectively.

Let us briefly describe $\gnbc$.
For all $S\in\inii$, consider the truncation by $S$:
$$
\pi_S\colon \rlx\lra\rlx,\qquad x=(x_i)_{i\in I}\longmapsto x_S=(y_i)_{i\in I},
$$
where $y_i=x_i$ for all $i\in S$ and $y_i=0$ otherwise.
Note that $\pi_S^{-1}(\be_S)=\be+H_S$, where $H_S$ is the convex subgroup defined in (\ref{HSHahn}).
The set $\gnbc$ is stratified as:
$$
\gnbc=\bigsqcup\nolimits_{S\in\inii}\op{nbc}(S),
$$
where \ $\op{nbc}(S)=\left\{x\in\pi_S\left(\rlx\right)\setminus \gq\mid x_T\in\gq \mbox{ for all } T\in\inii,\ T\subsetneq S \right\}$.

%The $\op{sme}$-class of these elements may be computed as follows:$$[\be]_\sme=\pi_S^{-1}(\be_S)=\be+H_S,$$where $H_S$ is the convex subgroup defined in (\ref{HSHahn}).

Now, let us describe $\gbc$.
%For all $S\in\inii$ consider unit vectors$$e_S=(e_j)_{j\in \I}\in\R^{I_S}_{\lx}\subset \rii,\qquad e_j=0\ \mbox{ for all }j\ne i_S,\quad e_{i_S}=1.$$
For all $b=(b_i)_{i\in I}\in\gq$, $S\in\inii$, denote
$$
b_S^-=((b_j)_{j\in S}\mid-1\mid0\cdots0)\in\R^{I_S}_{\lx},\qquad b_S^+=((b_j)_{j\in S}\mid1\mid0\cdots0)\in\R^{I_S}_{\lx},
$$
where $\pm1$ is placed at the $i_S$-th coordinate.
Then, $\gbc$ is constructed as:
$$
\gbc=\bigsqcup\nolimits_{S\in\inii}\left\{b_S^-,\,b_S^+\mid b\in\gq\right\}.
$$

The elements determined by the initial segment $S=\emptyset$ deserve a special notation:
$$
-\infty=(-1\mid0\cdots 0)=\min(\gsme),\qquad \infty^-=(1\mid0\cdots 0)=\max(\gsme),
$$
where $\pm1$ is placed at the $i_\emptyset$-th coordinate; that is, the first coordinate of $I_\emptyset=\{i_\emptyset\}+I$.
The notation for $\infty^-$ is motivated by the fact that this element is the immediate predecessor of $\infty$ in the set $\gsme\infty$. %Note that$$-\infty,\qquad \infty^-=\max(\gsme).$$

\section{Construction of the valuative tree}\label{secConstruct}

We keep with the notation of the previous section
$$\rll\subset\La:=\rii, \qquad \tz\subset\ttt:=\tla,$$
and we assume that the valuation $v$ takes values in a subgroup of  $\La$. %We recall that the inclusion $\g\subset\gq\subset\rlx$ preserves the rank.

%$$v(K^*)=\g\subset\gq\subset\gsme\subset\rll\subset \La,$$

%Thus, all valuations $\mu\in\ttt=\tla$ restrict to the valuation $v$ on $K$.

Since $\gsme$ is complete, we may extend the singular value function $\sval$ to the leaves of $\tz$. For a finite leaf $\nu\in\lfin(\tz)$, we agree that $\sval(\nu)=\infty$, while for an infinite leaf $\nu\in\li(\tz)$ we define
$$
\sval(\nu)=\sup\left\{\sval(\rho)\mid \rho\in(-\infty,\nu)\right\}\in\gsme.
$$

\subsection{Equivalence classes of commensurable extensions}
Let $\ttt_v$ be the set of equivalence classes of valuations on $\kx$ whose restriction to $K$ is equivalent to $v$.
It is well-known how to describe the subset $\ttt_v^\com \subset\ttt_v$ of the equivalence classes which are commensurable over $[v]$. %some valuation on $K$ equivalent to $v$.

By Proposition \ref{riiUniverse}, any such valuation is equivalent to some commensurable node  $\mu\in\ttt$; that is,
 a node belonging to the subtree:
$\tq:=\tgq\subset\ttt$.

Finally, it is easy to classify the nodes of $\tq$ up to equivalence.
Two commensurable valuations $\mu,\mu'\in\tq$ are equivalent if and only if $\mu=\mu'$. Indeed, %for any commensurable extension $\g\hk \Delta$, there exists a unique (order-preserving) embedding of $\Delta$ into $\gq$ such that  the composition $\g\hk \Delta\hk\gq$ is the canonical embedding. In particular,
if two subgroups
$$\g\subset\Delta\subset\gq,\qquad \g\subset\Delta'\subset\gq,$$
admit an order-preserving isomorphism  $\iota\colon \Delta \ism\Delta'$ which is the identity on $\g$, then necessarily $\Delta=\Delta'$ and $\iota$ is the identity mapping.

Therefore, we have a natural bijective mapping
$$
\tq\lra\ttt_v^\com,\qquad \mu\longmapsto[\mu].
$$

Since all leaves of $\ttt$ are commensurable, they  are leaves of $\tq$ too.
%Conversely, by Theorem \ref{maximal} any leaf $\mu$ of $\tq$ has $\kpm=\emptyset$; hence, it is a leaf of $\ttt$ too, . 
Therefore, both trees have the same leaves. More precisely, with the notation of Section \ref{subsecMLV}, we have:
\begin{equation}\label{EqualLeaves}
%\ars{1.4}
%\begin{array}{ccc}
\lfin(\tq)=\lfin(\ttt),\quad
\lui(\tq)=\lui(\ttt),\quad
\lci(\tq)=\lci(\ttt).
%\end{array}
\end{equation}

By Lemmas \ref{equiv=lim}, \ref{numerable} and \ref{specialCont}, every leaf in $\lui(\tq)$ is the stable limit of a countable family of nodes in $\tq$ with unbounded degree, and every leaf in $\lci(\tq)$ is the stable limit of an essential continuous family of nodes in $\tq$.

\subsection{Description of the valuative tree}
Consider the subtree
$$\ts=\left\{\rho\in\tz\mid \sval(\rho)\in\gsme\right\}.$$

Since $\gq\subset\gsme$, we have $\tq\subset\ts\subset\tz\subset\ttt$. In particular, from (\ref{EqualLeaves}) we deduce
$$%\begin{equation}\label{EqualLeaves2}
%\ars{1.4}
%\begin{array}{ccccc}
\lfin(\tq)=\lfin(\ts),\quad
\lui(\tq)=\lui(\ts),\quad
\lci(\tq)=\lci(\ts).
%\end{array}
$$%\end{equation}

\begin{theorem}\label{main2}
The mapping $\mu\mapsto[\mu]$ induces a bijection between $\ts$ and  $\ttt_v$.
\end{theorem}

\begin{proof}
Let $\mu$ be a valuation on $\kx$ whose restriction to $K$ is equivalent to $v$.  By Proposition \ref{riiUniverse}, $\mu$ is equivalent to some valuation in $\tz$. Thus, we may suppose $\mu\in\tz$. If $\mu$ is commensurable, then $\mu\in\tq\subset\ts$, so that $[\mu]$ is the image of some node of $\ts$.

Suppose that $\mu$ is incommensurable. Then, it is the last node of a finite MLV chain
$$ \mu_0\ \stackrel{\phi_1,\ga_1}\lra\  \mu_1\ \lra\ \cdots\ \lra\ \mu_{r-1}\ \stackrel{\phi_{r},\ga_{r}}\lra\ \mu_{r}=\nu.$$
If $\mu=\mu_0=\om_{a,\ga}$ has depth zero, then Corollary \ref{motivation2} shows that $\mu$ is equivalent to $\om_{a,\be}\in\ts$, where $\be\in\gsme$ is the representative of the class $[\ga]_\sme$.

If $\mu$ has a positive depth and $\mu_{r-1}\to\mu$ is an ordinary augmentation, then Corollary \ref{motivation2} shows that $\mu$ is equivalent to $[\mu_{r-1};\,\phi_r,\be]\in\ts$, where $\be\in\gsme$ is the representative of the class $[\ga_r]_\sme$.

If $\mu$ has a positive depth and $\mu_{r-1}\to\mu$ is a limit augmentation, then $\mu=[\aa;\,\phi_r,\ga_r]$ for some essential continuous family in $\tq$ and Corollary \ref{motivation2} shows that $\mu$ is equivalent to $[\aa;\,\phi_r,\be]\in\ts$, where $\be\in\gsme$ is the representative of the class $[\ga_r]_\sme$.

This proves that the mapping $\mu\mapsto[\mu]$ is onto.

Finally, let us show that the mapping $\mu\mapsto [\mu]$ is injective.
Suppose that $\mu,\nu\in\ts$ are equivalent. Then $\sval(\mu)\sim_\sme\sval(\nu)$ by  Proposition \ref{motivation}.
Since $\mu$ and $\nu$ belong to $\ts$, we have $\sval(\mu),\,\sval(\nu)\in\gsme$, so that necessarily $\sval(\mu)=\sval(\nu)$.

Also,  Proposition \ref{motivation} shows that $\kpm=\kpn$ and
$\nu(a)=\mu(a)$ for all $a\in\kx$ of degree less than $\deg(\mu)$.
This implies $\mu=\nu$ by comparison of their action on $\phi$-expansions, for any common key polynomial $\phi$ of minimal degree, having in mind that $\mu(\phi)=\sval(\mu)=\sval(\nu)=\nu(\phi)$.
\end{proof}\e

This subtree $\ts\subset\ttt$ shares many of the properties of $\ttt$ discussed in Sections \ref{secTreeLa} and \ref{secCtDepth}. Let us explicitely quote some of them.\e

$\bullet$\quad For all $\mu\in \ts$, the nodes of a MLV chain of $\mu$, except for (eventually) $\mu$ itself, are commensurable.
Thus, these nodes belong to $\ts$ and the depth of  $\mu$ can be described solely in terms of $\ts$. \e

$\bullet$\quad If $\mu\in \ts$ is an inner node and $\phi\in\kpm$, then we may build up ordinary augmentations in $\ts$:
$$
\nu=[\mu;\,\phi,\ga]\in \ts,\quad \ga\in\gsme,\ \ga>\mu(\phi).
$$
For any such augmentation, the interval $(\mu,\nu]\subset \ts$ may be described as
$$
(\mu,\nu]=\left\{[\mu;\,\phi,\dta]\mid \dta\in\gsme,\ \mu(\phi)<\dta\le\ga\right\}.
$$

$\bullet$\quad In particular, Proposition \ref{td=td} holds in $\ts$ too. There is a canonical bijection between $\kpm/\!\smu$ and the set of tangent directions of $\mu$ in the tree $\ts$.\e

$\bullet$\quad Let $\aa=(\rho_i)_{i\in A}$ be an essential continuous family in $\ttt$, and  $\phi\in\kpi(\aa)$ a limit key polynomial. Then, we may build up limit augmentations in $\ts$:
$$
\nu=[\aa;\,\phi,\ga]\in \ts,\quad \ga\in\gsme,\ \ga>\rho_i(\phi)\ \mbox{ for all }i\in A.
$$
By Lemma \ref{specialCont}, we may assume that all $\rho_i$ are commensurable. Thus, we may think that these limit augmentations are constructed solely from objects in the tree $\ts$.

For any such augmentation, we may  describe the following interval in $\ts$:
$$
\bigcap\nolimits_{i\in A}(\rho_i,\nu]=\left\{[\aa;\,\phi,\dta]\mid \dta\in\gsme,\ \rho_i(\phi)<\dta\le\ga\ \mbox{ for all }i\in A\right\}.
$$

$\bullet$\quad Every two nodes $\mu,\,\nu\in\ts$ have a greatest common lower node $\mu\wedge\nu\in\ts$.

Indeed, as remarked after Proposition \ref{GCN}, if neither $\mu\le\nu$ nor $\mu\ge\nu$, the greatest common lower node $\mu\wedge\nu\in\ttt$ is commensurable; thus, it belongs to $\ts$.

\subsection{Paths of constant depth in $\ts$}\label{subsecConstTs}
The main difference between $\ts$ and $\ttt$ lies in the fact that the paths of constant depth in $\ts$ are ``compact", thanks to the completeness of $\gsme$.

\subsubsection{Inner depth-zero nodes}
With the notation in Section \ref{subsecDepth0}, the inner depth-zero nodes of $\ts$ are of the form $\om_{a,\ga}$ for $a\in K$ and $\ga\in\gsme$. By (\ref{balls}), we have
$$
\om_{a,-\infty}=\om_{b,-\infty}\le\om_{c,\ga} \quad\mbox{ for all}\quad a,b,c\in K,\ \ga\in \gsme.
$$
Let us denote by $\minf:=\om_{a,-\infty}$ this minimal depth-zero valuation, which is independent of $a$.
By Theorem \ref{main}, this node $\minf$ is an absolute minimal node of $\ts$. We say that $\minf$ is the \emph{root node} of $\ts$. As a valuation, it works as follows:
$$
\minf\colon\kx\lra \left(\Z\times\g\right)\infty,\qquad f\longmapsto \left(-\deg(f),v(\lc(f))\right),
$$
where $\lc(f)$ is the leading coefficient of a nonzero polynomial $f$.
All valuations $\mu$ on $\kx$ satisfying $\mu(x)<\gq$ are equivalent to $\minf$ \cite[Thm. 2.4]{RPO}.

Since $\minf$ is incommensurable, it has a unique tangent direction. Actually, $$\kp(\minf)=\{x-a\mid a\in K\}=[x]_{\minf}.$$

All inner depth-zero nodes in $\ts$ are obtained as a single ordinary augmentation of the root node $\minf$:
$$
\om_{a,\ga}=[\minf;\,x-a,\ga]\quad\mbox{ for all }a\in K,\ \ga\in\gsme,\ \ga>-\infty.
$$
%by comparing the action of both valuations on $(x-a)$-expansions.

%Therefore, the set of all inner depth-zero nodes in $\ts$ beyond $\minf$ coincides with the set of all nodes in  constant-depth paths beyond $\minf$, determined by ordinary augmentations. 

In particular, the set of all inner depth-zero nodes is:
$$
\left\{\minf\right\}\cup\bigcup\nolimits_{a\in K}\pset_{\minf}(x-a).
$$
%All these constant-depth paths are weak because all key polynomials in $\kp(\minf)$ have degree one.
For any key polynomial $x-a\in \kp(\minf)$, the constant-depth path $\pset_{\minf}(x-a)$ is parametrized by the interval $(-\infty,\infty]\subset\gsme\infty$:

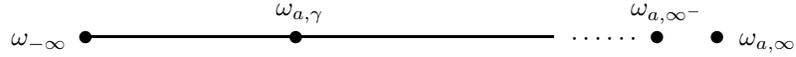
\begin{center}
\setlength{\unitlength}{4mm}
\begin{picture}(20,3.5)
\put(-1,1){$\bullet$}\put(-0.8,1.3){\line(1,0){15.6}}\put(18,1){$\bullet$}
\put(6,1){$\bullet$}\put(20,1){$\bullet$}
%\multiput(3,.5)(0,.25){22}{\vrule height2pt}
%\multiput(8,.9)(0,.25){9}{\vrule height2pt}
%\multiput(-.1,3)(.25,0){55}{\hbox to 2pt{\hrulefill }}
\put(-3.2,1){\begin{footnotesize}$\minf$\end{footnotesize}}
\put(15.4,1){\begin{footnotesize}$\cdots\cdots$\end{footnotesize}}
\put(17.4,2){\begin{footnotesize}$\om_{a,\infty^-}$\end{footnotesize}}
\put(21,1){\begin{footnotesize}$\om_{a,\infty}$\end{footnotesize}}
\put(5.6,2){\begin{footnotesize}$\om_{a,\ga}$\end{footnotesize}}
\end{picture}
\end{center}

Moreover, $\om_{a,\ga}$ is commensurable if and only if $\ga\in\gq\infty$, and it preserves the rank if and only if $\ga\in\gnbc\infty$.
The finite leaf $\om_{a,\infty}$ has an immediate predecessor node $\om_{a,\infty^-}$, represented by the valuation
$$
\om_{a,\infty^-}\colon\kx\lra \left(\Z\times\g\right)\infty,\qquad f\longmapsto \left(\ord_{x-a}(f),v(\init(f))\right),
$$
where $\init(f)$ is the first  nonzero coefficient of the $(x-a)$-expansion of $f\in\kx$.

The intersection of the depth-zero paths in $\ts$ determined by any two $a,b\in K$ may be computed as in Section \ref{subsecDepth0}:
$$
\pset_{\minf}(x-a)\cap\pset_{\minf}(x-b)=[\minf,\om_{a,v(b-a)}].
$$

\subsubsection{Ordinary augmentations}

Let $\mu\in\ts$ be an inner node and let $\phi\in\kpm$ be a key polynomial. The  constant-depth path $\pmph\subset\ts$, of all nodes in $\ts$ determined by an ordinary augmentation of  $\mu$ with respect to $\phi$, is parametrized by all $\ga\in\gsme\infty$ such that $\ga>\mu(\phi)$:

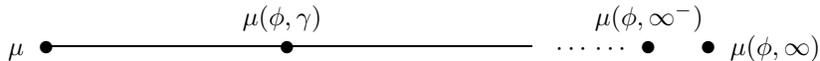
\begin{center}
\setlength{\unitlength}{4mm}
\begin{picture}(22,3.5)
\put(-2,1){$\bullet$}\put(-1.6,1.3){\line(1,0){16}}\put(18,1){$\bullet$}\put(20,1){$\bullet$}
\put(6,1){$\bullet$}
%\multiput(3,.5)(0,.25){22}{\vrule height2pt}
%\multiput(8,.9)(0,.25){9}{\vrule height2pt}
%\multiput(-.1,3)(.25,0){55}{\hbox to 2pt{\hrulefill }}
\put(-3,1){\begin{footnotesize}$\mu$\end{footnotesize}}
\put(16.5,1){\begin{footnotesize}$\cdots$\end{footnotesize}}
\put(15.2,1){\begin{footnotesize}$\cdots$\end{footnotesize}}
\put(16.5,2){\begin{footnotesize}$\mu(\phi,\infty^-)$\end{footnotesize}}
\put(21,1){\begin{footnotesize}$\mu(\phi,\infty)$\end{footnotesize}}
\put(4.8,2){\begin{footnotesize}$\mu(\phi,\ga)$\end{footnotesize}}
\end{picture}
\end{center}\e

Moreover, $\mu(\phi,\ga)$ is commensurable if and only if $\ga\in\gq\infty$, and it preserves the rank if and only if $\ga\in\gnbc\infty$.
The finite leaf $\mu(\phi,\infty)$ has an immediate predecessor node $\mu(\phi,\infty^-)$, represented by the valuation
$$
\mu(\phi,\infty^-)\colon\kx\lra \left(\Z\times\g\right)\infty,\qquad f\longmapsto \left(\ord_\phi(f),\mu(\init(f))\right),
$$
where $\init(f)$ is the first  nonzero coefficient of the $\phi$-expansion of $f\in\kx$.

The intersection of the constant-depth paths in $\ts$ determined by any two $\phi,\phi_*\in \kpm$ may be computed as in Section \ref{subsecConstDepthOrd}:
$$
\pmph\cap\pset_\mu(\phi_*)=\begin{cases}
\emptyset,&\quad\mbox{ if }\quad\phi\not\sim_\mu\phi_*,\\
\left(\mu,\mu(\phi,\ga_0)\right],&\quad\mbox{ if }\quad\ga_0=\mu(\phi-\phi_*)>\mu(\phi).
\end{cases}
$$

%Let us denote the set of all nodes in {\bf strong} constant-depth paths beyond $\mu$ by:$$\pmu\str=\bigcup\nolimits_{\phi\in\kpmz}\pmph\ \subset\  \ts,$$where $\kpmz$ is the subset of all strong key polynomials. We recall that $\mu\not\in\pmu\str$.

\subsubsection{Limit augmentations}

Finally, let $\aa=(\rho_i)_{i\in A}$ be an essential continuous family  in $\ts$, and let $\phi\in\kpi(\aa)$ be a limit key polynomial. Let $\mu\in\aa$ be the first valuation in the family. The completeness of $\gsme$ implies the existence of a minimal limit augmentation of $\aa$ in $\ts$ with respect to $\phi$; namely
$$
\mu_\aa:=[\aa;\,\phi,\ga_\aa],\qquad \ga_\aa:=\sup\left\{\rho_i(\phi)\mid i\in A\right\}\in\gsme.
$$
Note that $\ga_\aa>\rho_i(\phi)$ for all $i$, because $\aa$ has no maximal element. Also, $\ga_\aa<\infty$.

The following result follows immediately from Lemma \ref{allLKP}.

\begin{lemma}\label{muaa}
The value $\ga_\aa\in\gsme$ and the valuation $\mu_\aa\in\ts$ are independent of the choice of the limit key polynomial $\phi$ in $\kpi(\aa)$.
\end{lemma}

The constant-depth path $\paph\subset\ts$, of all nodes  determined by a limit augmentation of  $\aa$ with respect to $\phi$, is parametrized by all $\ga\in\gsme\infty$ such that $\ga\ge\ga_\aa$:

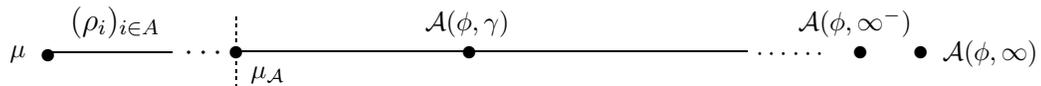
\begin{center}
\setlength{\unitlength}{4mm}
\begin{picture}(28,4)
\put(-2,0.9){$\bullet$}\put(4.25,1){$\bullet$}\put(25,1){$\bullet$}\put(27,1){$\bullet$}\put(12,1){$\bullet$}
\put(-1.6,1.25){\line(1,0){4}}\put(2.8,1){$\cdots$}\put(4.5,1.3){\line(1,0){17}}
\multiput(4.5,0.1)(0,.25){10}{\vrule height1pt}
\put(-1,2){$(\rho_i)_{i\in A}$}
%\multiput(8,.9)(0,.25){9}{\vrule height2pt}
%\multiput(-.1,3)(.25,0){55}{\hbox to 2pt{\hrulefill }}
\put(-3,1.1){\begin{footnotesize}$\mu$\end{footnotesize}}
\put(5,.4){\begin{footnotesize}$\mu_\aa$\end{footnotesize}}
\put(21.8,1){\begin{footnotesize}$\cdots$\end{footnotesize}}
\put(23,1){\begin{footnotesize}$\cdots$\end{footnotesize}}
\put(23.2,2){\begin{footnotesize}$\aa(\phi,\infty^-)$\end{footnotesize}}
\put(28,1){\begin{footnotesize}$\aa(\phi,\infty)$\end{footnotesize}}
\put(10.8,2){\begin{footnotesize}$\aa(\phi,\ga)$\end{footnotesize}}
\end{picture}
\end{center}\e

Note that $\mu_\aa\in\paph$. Again, $\aa(\phi,\ga)$ is commensurable if and only if $\ga\in\gq\infty$, and it preserves the rank if and only if $\ga\in\gnbc\infty$.
The finite leaf $\aa(\phi,\infty)$ has an immediate predecessor node $\aa(\phi,\infty^-)$, represented by the valuation
$$
\aa(\phi,\infty^-)\colon\kx\lra \left(\Z\times\g\right)\infty,\qquad f\longmapsto \left(\ord_\phi(f),\rha(\init(f))\right),
$$
where $\init(f)$ is the first  nonzero coefficient of the $\phi$-expansion of $f\in\kx$.

The intersection of the constant-depth paths in $\ts$ determined by any two $\phi,\phi_*\in \kpi(\aa)$ is an interval in $\ts$ which may be computed as in Section \ref{subsecConstDepthLim}:
$$
\paph\cap\pset_\aa(\phi_*)=[\mu_\aa,\aa(\phi,\rha(\phi-\phi_*)]\subset\ts.
$$

Since the set $\left\{\rho_i(\phi)\mid i\in A\right\}\subset\gq$ contains no maximal element, its supremum $\ga_\aa$ in $\gsme$ is incommensurable. Indeed, if $\ga_\aa\in\gq$, then it would admit an immediate predecessor  $\ga_\aa^-<\ga_\aa$, defined as $\ga_\aa^-=b_S^-$, for $b=\ga_\aa$ and $S=I$ (cf. Section \ref{subsecQcuts}). Since $\ga_\aa^-\in\gnbc$ is incommensurable,  it is still  an upper bound of the set $\left\{\rho_i(\phi)\mid i\in A\right\}$. This contradicts the minimality of $\ga_\aa$ as an upper bound of this set.

Thus, $\mu_\aa$ is incommensurable. In particular, it has a unique tangent direction.

Since $\ga_\aa<\infty$, Proposition \ref{extensionlim} shows that
$\phi$ is a key polyomial for $\mu_\aa$ of minimal degree. Actually, by \cite[Thm. 4.2]{KP} and Lemma \ref{allLKP}, we have
$$
\kp(\mu_\aa)=[\phi]_{\mu_\aa}= \left\{\phi+a\mid a\in\kx,\ \deg(a)<\mi, \ \rho_\aa(a)>\ga_\aa \right\}=\kpi(\aa).
$$

Also, all limit augmentations $\aa(\phi,\ga)$ are ordinary augmentations of $\mu_\aa$:
$$
\aa(\phi,\ga)=[\aa;\,\phi,\ga]=[\mu_\aa;\,\phi,\ga]\quad \mbox{ for all }\ga\in\gsme,\ \ga>\ga_\aa,
$$
by comparing the action of both valuations on $\phi$-expansions. Indeed, for all polynomials $a\in\kx$ of degree less than $\mi=\deg(\phi)$, we have $\mu_\aa(a)=\rha(a)$, by the definition of a limit augmentation.

The above picture might suggest that the interval $(\mu,\mu_A)$ is contained in a single constant-depth path beyond $\mu$. This is not the case.

By Lemma \ref{specialCont}, we may suppose that $\aa=\left(\rho_i\right)_{i\in A}$, with $\rho_i=[\mu;\,\chi_i,\be_i]$. Then, for each $i\in A$, the interval $(\mu,\rho_i]$ is contained in $\pset_\mu(\chi_i)$; however, for  $j>i$,  the valuation $\rho_j$ belongs to $\pset_\mu(\chi_j)$, but it does not belong to $\pset_\mu(\chi_i)$. Therefore, a more appropriate picture of this interval would be the following one:

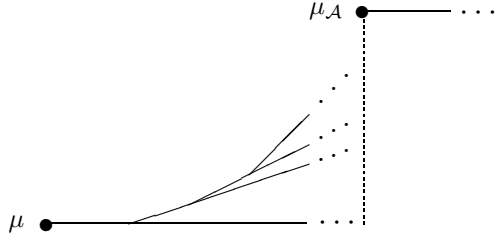
\begin{center}
\setlength{\unitlength}{4mm}
\begin{picture}(20,9.5)
\put(-2,0.9){$\bullet$}\put(-1.6,1.25){\line(1,0){8.5}}\put(7.3,1){$\cdots$}
\put(8.5,8){$\bullet$}\put(8.7,8.3){\line(1,0){3}}\put(12,8){$\cdots$}
\multiput(8.8,1.2)(0,.21){33}{\vrule height1pt}
%\multiput(8,.9)(0,.25){9}{\vrule height2pt}
%\multiput(-.1,3)(.25,0){55}{\hbox to 2pt{\hrulefill }}
\put(1,1.2){\line(3,1){6}}
\put(7.4,2.7){$\dot{}$}\put(7.82,2.85){$\dot{}$}\put(8.24,3){$\dot{}$}
\put(3,1.85){\line(2,1){4}}
\put(7.4,3.45){$\dot{}$}\put(7.82,3.66){$\dot{}$}\put(8.24,3.87){$\dot{}$}
\put(5,2.85){\line(1,1){2}}
\put(7.4,4.64){$\dot{}$}\put(7.82,5.06){$\dot{}$}\put(8.24,5.49){$\dot{}$}
%\put(18.9,5.26){$\dot{}$}\put(19.23,5.37){$\dot{}$}\put(19.56,5.5){$\dot{}$}

\put(-3,1.1){\begin{footnotesize}$\mu$\end{footnotesize}}
\put(7,8.2){\begin{footnotesize}$\mu_\aa$\end{footnotesize}}
\end{picture}
\end{center}

\subsection{Primitive nodes}

The constant-depth paths beyond a limit augmentation have completely analogous properties as the depth-zero paths. For the ease of the reader we include the depth-zero paths as a special case of the limit augmentation paths.\e

\nn{\bf Convention. }We admit the empty set $\aa=\emptyset$ as an essential continuous family in $\ts$. We agree that this family has $\ga_\aa=-\infty$, $\mu_\aa=\minf$, and
$$
\kpi(\aa)=\kp(\mu_\aa)=\left\{x-a\mid a\in K\right\},\quad
\pset_{\aa}(x-a)=\{\minf\}\,\cup\,\pset_{\minf}(x-a).
$$%\vskip0.2cm

\defn A \emph{primitive-limit} node in $\ts$ is the inner limit node $\mu_\aa$ associated to an essential continuous family $\aa$ in $\ts$. The set of primitive-limit nodes is in bijection with the set of equivalence classes of essential continuous families.

A \emph{primitive-ordinary\,} node in $\ts$ is an inner node $\mu\in\ts$ admitting strong constant-depth paths (cf. Section \ref{subsecConstDepthOrd}). That is, $\kpmz\ne\emptyset$, where
$$
\kpmz:=\{\phi\in\kpm\mid\deg(\phi)>\deg(\mu)\}.
$$
Since $\mu$ has key polynomials of different degrees, it is necessarily commensurable.

A \emph{primitive} node in $\ts$ is a node which is either primitive-limit or primitive-ordinary. Let us denote by  $\prim(\ts)$ the set of all primitive nodes. \e

%as $$\prim(\ts)=\primo(\ts)\,\sqcup\, \priml(\ts),$$where we distinguish the subsets of primitive-ordinary and primitive-limit nodes.\bs

By our convention, the root node $\minf$ is a primitive-limit node.

By Theorems \ref{main} and \cite[Thm. 4.7]{MLV}, the primitive-limit nodes cannot be obtained as an ordinary augmentation of a lower node.\e

\defn Let $\rho\in \ts$ be a primitive node. Then, we define
$$
\prh=\begin{cases}
\bigcup_{\phi\in\kp_{\op{str}}(\rho)}\pset_\rho(\phi),&\mbox{ if $\rho$ is primitive-ordinary},\\
\bigcup_{\phi\in\kpi(\aa)}\pset_\aa(\phi),&\mbox{ if $\rho=\mu_\aa$ is primitive-limit}.
\end{cases}
$$

We emphasize  that $\rho\in\prh$ if $\rho$ is a primitive-limit node, but $\rho\not\in\prh$ if $\rho$ is primitive-ordinary. However, in both cases, the arguments in Section \ref{subsecConstTs} show that
\begin{equation}\label{remind}
\mu\in\prh,\ \rho<\mu\ \imp\  \mu=[\rho; \phi,\sval(\mu)],
\end{equation}
for some $\phi\in\kpm$. If $\rho$ is primitive-ordinary, then necessarily $\phi\in\kpmz$.

\begin{theorem}\label{previous}
Let $\nu\in\ts$ be either an inner node, or a finite leaf. There exists a unique primitive node $\rho\in\prim(\ts)$ such that
$\nu\in\prh$. In other words,
$$
\ts\setminus\li(\ts)=\bigsqcup\nolimits_{\rho\in\prim(\ts)}\prh.
$$
\end{theorem}

\begin{proof}
If $\nu$ has depth zero, then $\nu$ belongs to $\pmi$, as we saw in Section \ref{subsecConstTs}.

If $\nu$ has a positive depth, then it is the last node of a finite MLV chain
$$ \mu_0\ \stackrel{\phi_1,\ga_1}\lra\  \mu_1\ \lra\ \cdots\ \lra\ \mu_{r-1}\ \stackrel{\phi_{r},\ga_{r}}\lra\ \mu_{r}=\nu.$$

If the last augmentation step $\mu_{r-1}\to\nu$ is ordinary, then $\mu_{r-1}$ is a primitive-ordinary node and $\nu\in\pset(\mu_{r-1})$. Indeed,
$\nu=[\mu_{r-1};\,\phi_r,\ga_r]\in\pset_{\mu_{r-1}}(\phi)$ and $\deg(\phi_r)>\deg(\mu_{r-1})$ by the definition of a MLV chain.

If $\mu_{r-1}\to\nu$ is a limit augmentation, then $\nu=[\aa;\,\phi_r,\ga_r]\in\pma$.

Therefore, the union of all sets $\prh$, for $\rho$ running on all the primitive nodes in $\ts$, covers $\ts\setminus\li(\ts)$. It remains only to show that
$$
\rho,\eta\in\prim(\ts), \ \ \rho\ne\eta\ \imp\ \prh\cap \pset(\eta)=\emptyset.
$$

Since $\ts$ is a tree, this is obvious if $\rho\not\le\eta$ and $\eta\not\le\rho$.

Suppose that $\rho<\eta$ and there exists $\mu\in\prh\cap\pset(\eta)$. By (\ref{remind}), the valuation $\mu\in\prh$ may be obtained after a single ordinary augmentation step: $\mu=[\rho;\,\phi,\sval(\mu)]$, for a certain $\phi\in\kpr$. Since $\eta$ belongs to the interval $(\rho,\mu)$, \cite[Lem. 2.7]{MLV} shows that $\eta=[\rho;\,\phi,\sval(\eta)]$ too. By Lemma \ref{propertiesAug}, $\deg(\eta)=\deg(\phi)=\deg(\mu)$.

This leads to a contradiction. Indeed, $\eta$ cannot be a primitive-limit node
because it is an ordinary augmentation of a lower node. Hence, $\pset(\eta)$ is the union of strong constant-depth paths and this implies $\deg(\eta)<\deg(\mu)$.
\end{proof}

\subsection{Stratification of $\ts$ by limit-depth}\label{subsecStrat}

Let $\rho\in\ts$ be a primitive-limit node. The \emph{inductive tree} with root $\rho$ is the subset
$\tind(\rho)\subset\ts$
formed by all inner nodes, or finite leaves in $\ts$, which may be obtained by a finite chain of {\bf ordinary} augmentations starting from $\rho$:
$$
\rho\ \stackrel{\phi_1,\ga_1}\lra\  \mu_1\ \lra\ \cdots
\ \lra\ \mu_{n-1}
\ \stackrel{\phi_{n},\ga_{n}}\lra\ \mu_{n}=\mu.
$$

We may consider the stratification by limit-depth
$$
\ts\setminus\li(\ts)=\bigsqcup\nolimits_{n\in\N_0}\ttt_n,
$$
where $\ttt_n$ is the subtree of all nodes in $\ts\setminus\li(\ts)$ whose limit-depth is equal to $n$.
These subtrees may be recursively constructed as:
$$
\ttt_0=\tind(\minf),\qquad \ttt_{n+1}=\bigsqcup\nolimits_{[\aa]\in\nni(\ttt_n)}\tind(\mu_\aa),
$$
where $\nni(\ttt_n)$ is the set of equivalence classes of essential continuous families in $\ttt_n$.

We could stratify $\li(\ts)$ in a similar way, but we must add a stratum corresponding to the infinite leaves with an infinite limit-depth. In \cite{ILD} we showed that such infinite leaves do exist.

\end{document}